\documentclass{amsart}
\usepackage[dvips]{graphicx}
\usepackage{amscd}
\usepackage{amsmath}
\usepackage{amsxtra}
\usepackage{amsfonts}
\usepackage{amssymb}

\newtheorem{theorem}{Theorem}[section]
\newtheorem{corollary}[theorem]{Corollary}
\newtheorem{lemma}[theorem]{Lemma}
\newtheorem{proposition}[theorem]{Proposition}

\theoremstyle{definition}
\newtheorem{definition}[theorem]{Definition}
\newtheorem{remark}[theorem]{Remark}

\newtheorem{example}[theorem]{Example}
\theoremstyle{remark}

\renewcommand{\theclaim}{\textup{\theclaim}}

\newtheorem*{acknowledgements}{Acknowledgements}

\numberwithin{equation}{section}

\def\openone%{\hbox{\upshape \small1\kern-3.3pt\normalsize1}}

{\mathchoice

{\hbox{\upshape \small1\kern-3.3pt\normalsize1}}

{\hbox{\upshape \small1\kern-3.3pt\normalsize1}}

{\hbox{\upshape \tiny1\kern-2.3pt\SMALL1}}

{\hbox{\upshape \Tiny1\kern-2pt\tiny1}}}

\makeatletter

\newbox\ipbox

\newcommand{\ip}[2]{\left\langle #1\, , \,#2\right\rangle}
\newcommand{\diracb}[1]{\left\langle #1\mathrel{\mathchoice

{\setbox\ipbox=\hbox{$\displaystyle \left\langle\mathstrut
#1\right.$}

\vrule height\ht\ipbox width0.25pt depth\dp\ipbox}

{\setbox\ipbox=\hbox{$\textstyle \left\langle\mathstrut
#1\right.$}

\vrule height\ht\ipbox width0.25pt depth\dp\ipbox}

{\setbox\ipbox=\hbox{$\scriptstyle \left\langle\mathstrut
#1\right.$}

\vrule height\ht\ipbox width0.25pt depth\dp\ipbox}

{\setbox\ipbox=\hbox{$\scriptscriptstyle \left\langle\mathstrut
#1\right.$}

\vrule height\ht\ipbox width0.25pt depth\dp\ipbox}

}\right. }

\newcommand{\dirack}[1]{\left. \mathrel{\mathchoice

{\setbox\ipbox=\hbox{$\displaystyle \left.\mathstrut
#1\right\rangle$}

\vrule height\ht\ipbox width0.25pt depth\dp\ipbox}

{\setbox\ipbox=\hbox{$\textstyle \left.\mathstrut
#1\right\rangle$}

\vrule height\ht\ipbox width0.25pt depth\dp\ipbox}

{\setbox\ipbox=\hbox{$\scriptstyle \left.\mathstrut
#1\right\rangle$}

\vrule height\ht\ipbox width0.25pt depth\dp\ipbox}

{\setbox\ipbox=\hbox{$\scriptscriptstyle \left.\mathstrut
#1\right\rangle$}

\vrule height\ht\ipbox width0.25pt depth\dp\ipbox}

} #1\right\rangle}

\newcommand{\cj}[1]{\overline{#1}}

\newcommand{\bz}{\mathbb{Z}}
\newcommand{\br}{\mathbb{R}}
\newcommand{\bc}{\mathbb{C}}
\newcommand{\bt}{\mathbb{T}}

\def\blfootnote{\xdef\@thefnmark{}\@footnotetext}

\renewcommand{\mod}{\operatorname{mod}}

\input xy
\xyoption{all}
\usepackage{amssymb}

\newcommand{\Hom}[2]{\text{Hom}(#1 #2)}
\newcommand{\End}[1]{\text{End}(#1)}

\def\R{\mathbb{R}}
\def\N{\mathbb{N}}

\def\T{\mathbb{T}}

\def\H{\mathcal{H}}

\def\-{^{-1}}

\def\C{\mathbb{C}}
\def\Z{\mathbb{Z}}

\begin{document}

\title[Harmonic functions]{The algebra of harmonic functions for a matrix-valued transfer operator}
\author{Dorin Ervin Dutkay}
\blfootnote{Research supported in part by a grant from the National Science Foundation DMS 0457491 and the Research Council of Norway, project number NFR 154077/420}
\address{[Dorin Ervin Dutkay]University of Central Florida\\
	Department of Mathematics\\
	4000 Central Florida Blvd.\\
	P.O. Box 161364\\
	Orlando, FL 32816-1364\\
U.S.A.\\} \email{ddutkay@mail.ucf.edu}
\author{Kjetil R\o ysland}
\address{[Kjetil R\o ysland]University of Oslo\\
Department of Mathematics\\
PO Box 1053, Blindern\\
NO-0316 Oslo\\
Norway}

\email{roysland@math.uio.no}
\thanks{
} 

\bibliographystyle{alpha}

\begin{abstract}
We analyze matrix-valued transfer operators. We prove that the fixed points of transfer operators form a finite dimensional $C^*$-algebra. For matrix weights satisfying a low-pass condition we identify the minimal projections in this algebra as correlations of scaling functions, i.e., limits of cascade algortihms. 
\end{abstract}
\maketitle \tableofcontents

\section{Introduction}\label{intro}
 In this paper we study a class of harmonic functions which come from the theory of dynamical systems. While there is an analogy to the classical theory of harmonic functions for Laplace operators, this analogy has two steps: First the continuous variable Laplace, or elliptic operators in PDE admit a variety of discrete approximations leading to random walk models, with variable coefficients PDEs corresponding to variable weight functions in the discrete version.

        As a result of work by numerous authors, and motivation from applications, there is now a rich harmonic analysis which is based on a certain transfer operator. It is based on path models, paths originating in a compact space, and with finite branching. The simplest instances come from endomorphisms in compact spaces $X$, onto, but with finitely branched inverses. A weight function on $X$ is prescribed, and the corresponding transfer operator (alias Ruelle operator) is denoted $R_W$.  

        Our motivation comes mainly from problems in wavelet analysis, but transfer operators are ubiquitous in applied mathematics: As it turns out, a study of transfer operators has now emerged as a subject of independent interest, with numerous applications. Examples: operator theory (locating the essential spectrum) \cite{CoIs91, Hel96}; $C^*$-algebras \cite{MuSo06}; groupoids \cite{MSS06, KuRe06}; infinite determinants \cite{BaRu96, Rue02}; anisotropic Sobolev spaces \cite{Bal05}; tiling spaces \cite{MaMu05}; number theory (zeta functions) \cite{Hen02, HiMa02, HiMa04, HMM05}; ergodic theory \cite{CHR97, DPU96}, dynamics \cite{CoRa98, CoRa03, Pol01}; optimization \cite{Hag05, KeLi99}; and quantum statistical mechanics \cite{Ara80, Rue92}.

       Our present paper uses tools from at least three areas, $C^*$-algebras, dynamics, and wavelets:

1.   We will be using both $C^*$-algebras, and their representations, especially an important family of matrix algebras.

2. Dynamical systems. The connection between the algebras are relevant to the kind of dynamical systems which are built on branching-laws. The reason for this is that the spectral properties of the transfer operator and of the associated eigenspace algebras are connected to the ergodic properties of the dynamical system.

3. Wavelet analysis. The connection to basis constructions using wavelets is this: The context for wavelets is a Hilbert space $H$, where $H$ may be $L^2(\br^d)$ where $d$ is a dimension, $d=1$ for the line (signals), $d=2$ for the plane (images), etc. The more successful bases in Hilbert space are the orthonormal bases ONBs, but until the mid 1980s, there were no ONBs in $L^2(\br^d)$ which were entirely algorithmic and effective for computations.

      Originating with \cite{Law90}, a popular tool for deciding whether or not a candidate for a wavelet basis is in fact an ONB uses a wavelet version of the transfer operator, still based on $n$-fold branching laws, but now with the branching corresponding to frequency bands. The wavelet Ruelle operator weights input over $n$ branching possibilities, and the weighting is assigned by a prescribed scalar function $W$, the modulus squared of a low-pass filter function, often called $m_0$. We are interested in the top part of the spectrum of $R_W$, a distinguished eigenspace for $R_W$; an infinite-dimensional version of the so called Perron-Frobenius problem from finite-dimensional matrix theory.

      This is especially useful for wavelets that are initialized by a single function, called the scaling function. These are called the multiresolution analysis (MRA) wavelets, or for short the MRA-wavelets. But there are multiwavelets (i.e., more than one scaling function) for example for localization in frequency domain, where the MRA-wavelets do not suffice, where it will by necessity include more than one scaling function. And then asking for an ONB is not feasible, but instead frame wavelets are natural.

       We attack this problem by introducing a matrix version of the weight function $W$; so our $W$ is no longer scalar valued, but rather matrix-valued; the size of the matrices depending on an optimal number of scaling functions. We show that the space of harmonic functions for the $W$-Ruelle operator with matrix valued weights acquires the structure of a $C^*$-algebra.  It serves to decide ONB vs frame properties, and the stability properties needed in applications, for example in the analogue to digital signal problem.

A transfer operator, also called Ruelle operator, is associated to a finite-to-one endomorphism $r:X\rightarrow X$ and a  weight function $W:X\rightarrow[0,\infty)$, and it is defined by
$$R_Wf(x)=\sum_{r(y)=x}W(y)f(y)$$
for functions $f$ on $X$. Here $X$ is a compact Hausdorff space, and $W$ is a non-negative continuous function on $X$. The function $W$ is said to be {\it normalized} if
$$\sum_{y\in r^{-1}(x)} W(y)=1,\quad(x\in X).$$ 
\par
Transfer operators have been extensively used in the analysis of discrete dynamical systems \cite{Baladi} and in wavelet theory \cite{BraJo02}.
\par
In multivariate wavelet theory (see for example \cite{Jiang} for details) one has an expansive $n\times n$ integer matrix $A$, i.e., all eigenvalues $\lambda$ have $|\lambda|>1$, and a multiresolution structure on $L^2(\br^n)$, i.e., a sequence of closed subspaces $\{V_j\}_{j\in\bz}$ of $L^2(\br^n)$ such that
\begin{enumerate}
\item $V_j\subset V_{j+1}$, for all $j$;
\item $\cup_j V_j$ is dense in $L^2(\br^n)$;
\item $\cap_j V_j=\{0\}$;
\item $f\in V_j$ if and only if $f( (A^T)^{-1}\cdot)\in V_{j-1}$;
\item There exist $\varphi_1,...,\varphi_d\in V_0$ such that $\{\varphi_k(\cdot-j)\,|\,k\in\{1,...,d\},j\in\bz^n\}$ forms an orthonormal basis for $V_0$.
\end{enumerate}

The functions $\varphi_1,...,\varphi_d$ are called {\it scaling functions}, and their Fourier transforms satisfy the following {\it scaling equation}:
$$\hat\varphi_i(x)=\sum_{j=1}^dm_{ji}(A^{-1}x)\hat\varphi_j(A^{-1}x),\quad(x\in\br^n,i\in\{1,...,d\}),$$
where $m_{ji}$ are some $\bz^n$-periodic functions on $\br^n$. 
\par
The orthogonality of the translates of $\varphi_i$ imply the following {\it QMF equation}:
$$\sum_{Ay=x\mod\bz^d}m^*(y)m(y)=1,\quad(x\in\br^n/\bz^n),$$
where $m$ is the $d\times d$ matrix $(m_{ij})_{i,j=1}^d$. 
When the translates of the scaling functions are not necessarily orthogonal one still obtains the following relation: if we let
$$h_{ij}(x):=\sum_{k\in\bz^n}\cj{\hat\varphi_i}(x+k)\hat\varphi_j(x+k),\quad(x\in\br^n/\bz^n),$$
then the matrix $h=(h_{ij})_{i,j=1}^d$ satisfies the following property:
\begin{equation}\label{eqrh}
Rh(x):=\sum_{Ay=x\mod\bz^n}m^*(y)h(y)m(y)=h(x),\quad(x\in\br^n/\bz^n),
\end{equation}
i.e., $h$ is a fixed point for the {\it matrix-valued transfer operator
} $R$. The fixed points of a transfer operator are also called {\it harmonic functions} for this operator. Thus the orthogonality properties of the scaling functions are directly related to the spectral properties of the transfer operator $R$. 
\par
This motivates our study of the harmonic functions for a matrix-valued transfer operator. The one-dimensional case (numbers instead of matrices) was studied in \cite{BraJo02,Dugale,Duspec}. These results were then extended in \cite{DuJoifs,DuJomart}, by replacing the map $x\mapsto Ax\mod\bz^n$ on the torus $\bt^n:=\br^n/\bz^n$, by some expansive endomorphism $r$ on a metric space.
\par
Here we are interested in the case when the weights defining the transfer operator are matrices, just as in equation \eqref{eqrh}. We keep a higher level of generality because of possible applications outside wavelet theory, in areas such as dynamical systems or fractals (see \cite{DuJoifs,DuJomart}). However, for clarity, the reader should always have the main example in mind, where $r:x\mapsto Ax\mod\bz^n$ on the $n$-torus. The quotient map $\br^n \rightarrow
\br^n / \bz^n$  defines a simply connected covering space,
and $r$ lifts to the dilation $\tilde r(x)=Ax$ on $\br^n$.
\par
In Section \ref{tran} we introduce the main notions. Since we are especially interested in continuous harmonic maps we used the language of vector bundles (see also \cite{PackerRieffel1}). Packer and Rieffel introduced projective multiresolution analyses (PMRA's) in
\cite{PackerRieffel1}. In their formalism, the scaling spaces correspond
to sections in vector bundles over $\bt^n$.  Motivated by their work, we
introduce transfer operators that act on bundlemaps on vectorbundles. This
gives a way to construct new PMRA's.\par
 In Section \ref{peri} we perform a spectral analysis of the transfer operator, prove that it is quasi-compact and give an estimate of the essential spectral radius. Section \ref{c*al} contains one of the main results of the paper: the continuous harmonic functions form a $C^*$-algebra, with the usual addition, multiplication by scalars and adjoint, and with the multiplication defined by a certain spectral projection of $R$. 
\par
In Section \ref{refi} we define the refinement operator. This operates at the level of the covering space $\tilde X$, and extends the usual refinement operator from wavelet theory (see \cite{BraJo02}). The correlations of the scaling functions correspond to fixed points of the refinement operator. We give the intertwining relation between the transfer operator and the refinement operator in Theorem \ref{thinter}. With this relation we will see that the fixed points of the transfer operator correspond to the fixed points of the refinement operator, i.e., to scaling functions.
\par
In Section \ref{exam} we consider the case of {\it low-pass} filters. In the one-dimensional case the low-pass condition amounts to $m(0)=1$, i.e., $m(0)$ is maximal. This is needed in order to obtain solutions of the scaling equation in $L^2(\br^n)$. In the matrix case, the appropriate low-pass condition was introduced in \cite{Jiang} under the name of the {\it  $E(l)$-condition} (see Definition \ref{defel}). If this condition is satisfied we show that the associated scaling functions exist in our case as pointwise limits of iterates of the refinement operator (Theorem \ref{prop5-1}). Their correlation functions give minimal projections in the algebra of continuous harmonic functions (Theorem \ref{thh-v}).
\par
In Theorem \ref{computation} we show that, if the peripheral spectrum of the transfer operator is ``simple'', and we make an appropriate choice of the starting point, then the iterates of the refinement operator converge strongly.

\section{Transfer operators} \label{tran}
We describe now our setting, starting with the main example, the one used in multivariate wavelet theory.
\begin{example}
Let $\tilde X = \R^n$ and let $G \subset \tilde X$ be a discrete subgroup
such that
$\tilde X /G$ is compact, i.e.
$G$ is a full-rank lattice, $G \simeq \Z^n$  and $X:=\tilde X  /G \simeq
\T^n$
and  let $A \in GL(\br^n)$ such that $A G \subset G$
and such that if $\lambda$ is an eigenvalue of $A$ then $|\lambda | > 1$. Then there exists a norm $\|\cdot\|_A$ on $\br^n$ and a number $\rho>1$ such that $\|Ax\|_A\geq \rho \|x\|_A$ for all $x\in\br^n$.
\par

Indeed, let $1  < \rho < \min |\lambda|$ and define $\| \cdot \|_A : \R^n
\rightarrow [ 0, \infty )$
as
$$
\| x\|_A = \sum_{j = 0}^\infty \rho^j \| A^{-j} x\|
$$
By the spectral radius formula and the root test,
this converges and defines a norm on $\R^n$ such that
$\| Ax\|_A \geq \rho \| x \|_A$ for all $x\in\R^n$.

Let $p : \tilde X \rightarrow X$ denote the quotient map and define  
a map $r : X \rightarrow X$ as 
$r(p(x)) = p(Ax)$. The map $r$ is a regular covering map with deck transformations $G/AG$.
We have that $\# G/AG = | \det(A)|$ and denote this number by $q:=\# G/AG$.
If $\mu$ is the normalized Haar measure on $X$,
then 
\begin{equation}\label{eqstrinv}
\int_X  f d \mu =   \int_X  q\- \sum_{rs=t} f(s) d\mu(t), 
\end{equation}
i.e., $\mu$ is {\it strongly $r$-invariant}.
Finally, there exists a normalization of the Haar measure on $\tilde X$, say $\tilde \mu$ 
such that 
$\int_{\tilde X} f d \tilde \mu = \int_X{\sum_{g \in G} f \circ g}d \mu$ for every 
$f \in C_c(\tilde X)$.
\end{example}

  We are mainly concerned with the above situation, but our results 
  apply to a more general setting: 
  \subsection{The setting}
  \par
    Let $\tilde X$ be a  locally compact metric space  with an  
  isometric and properly discontinuous 
  action of a discrete group 
  $G\times \tilde X \ni  g, \tilde x \mapsto g\tilde x$
  such that $X:=\tilde X / G$  is  compact. 
  
  Let $p:\tilde X\rightarrow X$ be the quotient covering map.
  \par
  Suppose $\tilde r$
  is  a  strictly expansive homeomorphism on $\tilde X$, i.e., for some $0<\theta<1$ 
  $$d(\tilde rx,\tilde ry)\theta\geq  d(x,y),\quad(x,y\in\tilde X).$$

  Let $\tilde x_0$ be the fixed point of $\tilde r$ and $x_0:=p(\tilde x_0)$.
  
  Assume that there exists
  $A \in \text{End} ( G )$ such that  
  $AG \subset G$ is a normal subgroup, and 
  $\tilde rg\tilde x = (Ag) \tilde r \tilde x$  for every 
  $g \in G$, $\tilde x \in \tilde X$. 
  
  Let $r:X\rightarrow X$, $r(p(x))=p(r(x))$, for all $x\in\tilde X$.
  \par
  
  Moreover, let $\tilde \mu$ and $\mu$ be  regular 
  measures on $\tilde X$ and $X$
  such that $  \mu$ is strongly $r$-invariant as in \eqref{eqstrinv}, and  
  \begin{equation}\label{eqpermu}
\int_{\tilde X} f d \tilde \mu = \int_X{\sum_{g \in G} f \circ g}d \mu,\quad (f \in C_c(\tilde X)).
\end{equation}
\par
Let $\rho : \xi \rightarrow X$ be  a Lipschitz continuous  $d$-dimensional 
complex vectorbundle over $X$  with 
a hermitian metric,
 i.e. a  map 
$\langle \cdot , \cdot \rangle : \xi \times \xi \rightarrow \C$ that restricts to 
positive definite sesquilinear forms on each fiber. 
We say that $\xi$ is Lipschitz continuous if there exists 
a trivializing system of bundlemaps 
$\phi_{U_i} : \xi |_{U_i} \rightarrow U_i \times \C^d $
such that every $\phi_{U_i} \phi_{U_j}\- $ is Lipschitz continuous 
on $(U_i \cap U_j) \times \C^d$.
\par
Let $S$ denote the continuous sections in $\xi$. By the Serre-Swan theorem 
\cite{Atiyah1}, this is a 
projective $C(X)$-module and the endomorphisms $F$ on $S$ are exactly the
bundlemaps on $\xi$ acting as $Fs(x) = F(x)s(x)$. In fact,  $\text{End}(\xi)$
equipped with the norm $\|F\|_\infty = \sup_{x \in X} \|F(x)\|$ and the pointwise involution with 
respect to the form $\langle \cdot , \cdot \rangle_x  : \xi_x \times \xi_x \rightarrow \C$
 is a $C^*$-algebra.

 Let $p^*\xi$ be the pull-back of the vector bundle $\xi$ by the map $p$, i.e., 
 $$p^*\xi:=\{(x,v)\in\tilde X\times \xi\,|\, p(x)=\rho(v)\}.$$
{\bf Assumption:} We assume that $p^* \xi$ is a trivial $d$-bundle. \par
This is always the case if $\tilde X$ is contractible. We claim that we can take $p^*\xi$ to be {\it the} trivial bundle $\tilde X\times\C^d$ with the canonical Hermitian inner product.
\par
Given $d$-linearly independent sections $s_1, \dots ,s_d$ in $p^*\xi$ 
such that $ \phi_{U_i} s_i|_{U_i}$ is Lipschitz continuous for every $i,j$
we get a Lipschitz continuous bundle isomorphism 
$\sigma :   \tilde X \times \C^d \rightarrow  p^* \xi$ such that
$\sigma e_i = s_i$, where $e_i$, $1\leq i\leq d$ are the canonical sections in the trivial bundle $\tilde X\times\C^d$. 
We equip   the product bundle with the standard inner product on $\C^d$
and $p^*\xi$ with pullback of the inner product on $\xi$.  
Let $u_x | \sigma_x | = \sigma_x$ denote the pointwise polar decomposition. 

 $\sigma^*\sigma$ is positive and invertible, i.e. the holomorphic 
functional calculus on Banach algebras  tells us that we can apply the square root 
and still get a bounded operator. This means that 
$|\sigma|$ is a bounded operator on the Lipschitz continuous $\tilde X \rightarrow \C^d$.
Now $x,v \mapsto  u_xv$ defines a Lipschitz continuous bundle isomorphism 
$\tilde X \times \C^d \rightarrow p^*\xi$ that  is isometric over each $x$.  
This means that we can identify $p^*\xi$ 
with the trivial product bundle equipped with the ordinary inner product.

 Let $S_1 \subset S$ denote the set of 
  $s \in S$ such that $p^*s \in C_b(\tilde X, \C^d)$ 
 is Lipschitz continuous. $S_1$ is a projective $\text{Lip}_1(X)$
  module of finite rank and a Banach space  with the norm 
  $$
  \|s\|_1 = \sup_{x \in X} \| s(x)\| + 
  \sup_{x,y\in \tilde X}\frac{\| p^*s  (x) - p^*s(y) \|_{\C^d} }
  {d(x,y)}.
  $$
  Let $L$ denote the endomorphisms  on $S_1$. $L$ is dense in $\text{End}(\xi)$ and 
  a Banach 
  algebra  with the usual operator norm
  $$ \| g \|_1 = \sup_{ \{ s \in S_1 | ~ \| s \|_1 \leq 1\}} \| gs \|_1.$$

Let $r^*\xi$ denote the pullback of $\xi$ along $r$ (see \cite[1.1]{Atiyah1}) and 
let $$m : r^*\xi \rightarrow \xi$$ be a bundle map.
The pullback  of the inner product on $\xi$ gives an inner product on $r^*\xi$ 
and we get a unique $m^* : \xi \rightarrow r^*\xi$ 
such that $\langle v , m w \rangle  = \langle m^* v, w \rangle$ for every
pair $v \in \xi_x$ $w \in r^* \xi|x$.

Let $g_1 , \dots, g_q $ be a complete system of representatives for the right cosets of $AG$ in $G$, with $g_1=1$, and define 
$\psi_1, \dots , \psi_q : \tilde X \rightarrow \tilde X$ as 
\begin{equation}\label{eqpsi}
\psi_i =  \tilde r\- \circ g_i,
\end{equation}
we see that $d( \psi_i( x) ,\psi_i( y)) \leq \theta d(x,y)$, and $\{p\psi_1x,...,p\psi_qx\}=r^{-1}(px)$ for all $x\in\tilde X$.
\begin{definition}
We define the {\it transfer operator} associated to $m$  as the operator
$R$ acting on $h\in\text{End}(\xi)$ such that 
$$
(Rh) (px) =   \sum_{ i = 1}^q    m^* (p\psi_i x) h (p \psi_i x)   m (p \psi_i x),\quad(x\in\tilde X).
$$ 
An element $h\in\End\xi$ is called {\it harmonic} for the transfer operator $R$, if $Rh=h$.
Define 
$$
m^{(k)}(x) := m(px) m(rpx) \dots    m (r^{k-1} px).
$$
Moreover, let $\Omega_k := \times_{i= 1}^k \{1 , \dots , q\} $, $\omega \in \Omega_k$ and let 
$\omega_j$ denote the $j$'th coordinate of $\omega$.
\end{definition}
\par
A computation yields the following identity: 
$$ (R^kh)(px)\\  = 
\sum_{ \omega \in \Omega_k}    
m^{(k)*}  (\psi_{\omega_k} \dots \psi_{\omega_1} x) 
h(  p\psi_{\omega_k} \dots \psi_{\omega_1} x     )   
m ^{(k)}( \psi_{\omega_k} \dots \psi_{\omega_1} x ).
$$

We will use repeatedly the following Cauchy-Schwarz type inequality:
\begin{lemma}\label{CauchySchwarz}
  If  $m_1 , \dots,  m_k ,\tilde m_1 , \dots,  
  \tilde m_k ,h_1,  \dots,  h_d \in M_d(\C)$ 
  then  $$\| \sum_{i=1}^k \tilde m_i^* h_i m_i \| \leq 
  \| \sum_{i=1}^k \tilde m_i^* \tilde m_i \|^{1/2} ~   \| \sum_{i=1}^k  m_i^*  m_i \|^{1/2}
  \max_{1\leq j\leq k} \| h_j\|.$$
  \end{lemma} 
\begin{proof}
 The space $E := \oplus_{j=1}^k M_d(\C)$ with the sesquilinear map 
  $$
  \langle \cdot , \cdot \rangle : E \times E \rightarrow M_d(\C),
  \quad
  \langle a_1\oplus \dots \oplus a_k , b_1 \oplus \dots \oplus b_k \rangle = \sum_{i=1}^k a_i^* b_i
  $$
  form an $M_d(\C)$-Hilbert module. 

  If $H \in \text{End}_{M_d(\C)}(E)$ then by the Cauchy-Schwarz inequality for Hilbert Modules 
  $$
  \|\langle \eta, H \zeta \rangle\| \leq 
  \|  H \zeta\| ~ \|\langle \eta, \eta \rangle \|^{1/2} \leq \| H\| ~ 
  \| \langle \zeta, \zeta \rangle \|^{1/2}  ~   \| \langle \eta, \eta \rangle \|^{1/2}. 
  $$
  The claim  follows with  $H (a_1\oplus\dots\oplus a_k)=(h_1a_1\oplus\dots\oplus h_ka_k)$,  $\zeta = m_1 \oplus \dots \oplus m_k$ 
  and    
  $\eta = \tilde m_1 \oplus \dots \oplus \tilde m_k$.
\end{proof}

If there exists $h \in \End\xi$ that is positive, invertible and harmonic with respect to $R$, then there exists a $c > 0$
such that $c 1 \leq h$; this implies that for all $k\geq0$, $h=R^kh\geq cR^k1 \geq 0$, and with Lemma \ref{CauchySchwarz}, for all $h_0\in\End\xi$:   
\begin{align*}
  \| R^k h_0 \|_\infty & \leq \| R^k 1 \|_\infty \| h_0 \|_\infty \leq c\- \|h\|_\infty  \|h_0\|_\infty,  
\end{align*}
so the existence of such an element implies that $\sup_k \|R^k\|_\infty < \infty$.

We assume from now on that $m$ satifies the following conditions. 
\begin{enumerate}
\item $m : r^*\xi \rightarrow \xi$ is Lipschitz continuous.
\item $ \sup_k \| R^k \|_\infty < \infty$
\end{enumerate}
\begin{remark}\label{remgraph} Since $p^*\xi$ is the trivial bundle, we have
$p^*r^* \xi = \tilde r ^* p^* \xi =  \text{graph}(\tilde r) \times \C^d$. 
This means that $p^*m$ defines a bundle map  $\text{graph}( \tilde r)  \times \C^d \rightarrow \tilde X \times \C^d$. 
Now $p^*m(x,y,v) = (x, m_0(x)v)$ for some map $m_0 : \tilde X \rightarrow M_d(\C)$. $p^*m^*$ defines a map 
$\tilde X \times \C^d \rightarrow \text{graph}( \tilde r) \times \C^d $ such that 
$p^*m^* (x,v) = (x,\tilde rx, m_0^*(x) v)$ where $m^*_0(x)$ is the unique 
element map $\tilde X \rightarrow M_d(\C)$ such that $\langle m_0^*(x) v ,  w \rangle =
\langle v , m_0(x) v \rangle$ for every $x \in \tilde X$
and $v,w \in \C^d$. Therefore we can identify $p^*m$ with $m_0$.
\end{remark}

\section{The peripheral spectrum}\label{peri}
  The {\it essential spectral radius} of a bounded 
linear operator $T$ on a Banach space is 
the infimum of positive numbers  $\rho \geq 0$ such that
$\lambda \in sp(T)$ and 
$|\lambda| > \rho$ implies $\lambda$ is an isolated eigenvalue of 
finite multiplicity. 
If the essential spectral radius of $T$ is strictly less than the spectral radius 
then $T$ is said to be quasicompact. 

Whenever $Y$ is a complete metric space and $ Z \subset Y$ we define the {\it 
Ball measure of noncompactness }of  $Z$, say $\gamma_Y(Z)$  to be 
$$
\gamma_Y(Z) = \inf
\{ \delta > 0 |\mbox{ There exist } z_1, \dots , z_k \in Z \mbox{ s.t. }Z \subset \cup_i B(z_i,\delta ) \}.
$$
If $B$ is a  Banach space with unit ball $B_1$ and $R$  is a bounded  
 linear operator on $B$, define 
$$
\gamma(R) = \gamma(R(B_1)).
$$
The following theorem is due to Nussbaum \cite{Nussbaum}
\begin{theorem}
$\lim_k \gamma(R^k)^{1/k}$ exists and equals the essential spectral radius of $R$. 
\end{theorem}
The corollary is due to Hennion \cite{Hennion}
\begin{corollary} \label{quasicompact}
  If $(L,\|\cdot \|)$ is a Banach space with a  second norm $| \cdot|$
and an  operator $R$ 
  such that

  \begin{enumerate}
    \item $R : (L, \| \cdot \|) \rightarrow 
      (L , | \cdot |)$ is a compact operator.   
    \item For every $n \in \N$ there exist positive numbers $r_n$ and $R_n$ such that 
      $\liminf_n r_n^{1/n} \leq r$ and 
      $$
      \| R^n  f \| \leq R_n |f| + r_n \|f\|,
      $$
      then the essential spectral radius of $R$ is less than $r$.
  \end{enumerate}
\begin{proof}
Let $B_1 = \{b \in L |~ \| b \| = 1\}$, 
$D(z,\delta) = \{y \in L | ~|y-z | < \delta   \}$
and 
$B(z,\delta) = \{y \in L | ~\|y-z \| < \delta   \}$.
Since $R(B_1)$ is relatively compact with respect to $(L, |\cdot|)$ 
there exists a sequence $b_1, \dots , b_m \in B_1$ for every $\delta > 0$ 
such that 
$$
R(B_1) \subset \cup^m_{i=1} D( R b_i, \delta) \cap B(0,\|R\|).
$$
If $a \in D( R b_i, \delta) \cap B(0,\|R\|)$ then 
$$
\| R^n Rb_i - R^n a \| \leq R_n |Rb_i - a  |  + r_n \|Rb_i - a\| \leq R_n \delta + 2 
\|R\| r_n. 
$$
Thus $R^{n+1}(B_1)$ can be covered by $\{B(R^{n+1} b_i , R_n \delta + 2 
\|R\| r_n   )\}_i$.
Such a sequence can be found for arbitrary $\delta > 0 $, i.e. 
$\gamma(R^{k+1}) \leq 2 \| R\| r_n$ and 
$$\lim_n \gamma(R^n)^{1/n} \leq \liminf_n (2\| R \| r_n    )^{1/n}  \leq   r.$$

\end{proof}
\end{corollary}
We intend to give an estimate of  the essential spectral radius of $R|_L$. First we need some lemmas.

\begin{lemma} \label{key estimate 1}
  If    $m \circ p $ is  Lipschitz continuous, there exists a $D > 0$ such that  
\begin{align*}
& ( \sum_{\omega \in \Omega_k}   \|  m^{(k)} ( \psi_{\omega_k}  \dots \psi_{ \omega1} x)   -  
       m^{(k)} ( \psi_{\omega_k}  \dots \psi_{ \omega1} y)   \|^2 )^{1/2}  \leq   d(x,y)D.
\end{align*}
For every $x,y \in \tilde X$ and $k \in \N$.
  \begin{proof}
  Since $\tilde r$ is expansive, the maps $\psi_j$ are contractive, with contraction constant $\theta$.
    $$
       \sum_{\omega \in \Omega_k}   \| m^{(k)} ( \psi_{\omega_k}  \dots \psi_{ \omega_1} x)  -  
      m^{(k)} ( \psi_{\omega_k}  \dots \psi_{ \omega_1} y) \|^2  $$
      $$
       =   \sum_{\omega \in \Omega_k}   \| \sum_{j = 1}^k 
         m^{(k-j)} ( \psi_{\omega_k } \dots \psi_{\omega_1} x) 
	    ( m(p \psi_{\omega_j} \dots \psi_{\omega_1} x ) -  
	 m(p \psi_{\omega_j} \dots \psi_{\omega_1} y ) )
         m^{(j-1)}( \psi_{\omega_{j-1}}  \dots  \psi_{\omega_1} y ) \|^2 $$
       
  %      $$ \leq & 
  %       \sum_{\omega \in \Omega^{(k)}} (  \sum_{j = 1}^k \| 
  %       m ( \psi_{\omega_k } \dots \psi_{\omega_1} x)  \dots m(\psi_{\omega_{j+1}} \dots \psi_{\omega_1} x)\|
  %        \\  &   ~\| m( \psi_{\omega_j} \dots \psi_{\omega_1} x ) 
  %        -  m( \psi_{\omega_j} \dots \psi_{\omega_1} y ) \|
  %       ~\|m( \psi_{\omega_{j-1}}  \dots  \psi_{\omega_1} y ) \dots m(\psi_{\omega_1} y) \|)^2 
  %       \\
       $$   \leq  
        \text{const}  \sum_{\omega \in \Omega_k} (  \sum_{j = 1}^k \| 
         m^{(k-j)} ( \psi_{\omega_k } \dots \psi_{\omega_1} x)\|
             \theta^j d(x,y)
         ~\|m^{(j-1)}( \psi_{\omega_{j-1}}  \dots  \psi_{\omega_1} y )  \|)^2.
         $$
    By the Cauchy-Schwarz inequality , this is dominated by 
    $$
       \text{const} ~d(x,y)^2   \sum_{\omega \in \Omega_k} ( (\sum_{i^=1}^k \theta^{2i})^{1/2}    (\sum_{j = 1}^k \| 
         m^{(k-j)} ( \psi_{\omega_k } \dots \psi_{\omega_1} x)\|^2 
          \|m^{(j-1)}( \psi_{\omega_{j-1}}  \dots  \psi_{\omega_1} y )  \|^2 )^{1/2} )^2.
    $$
   Since $ \| H \| \leq \text{tr} (H) \leq d \| H\|$ for any positive $d\times d$ matrix, we obtain 
    \begin{equation}\label{eqhi}\sum_i \|H_i\|^2=\sum_i \| H^*_i H_i \| \leq \sum_i \text{tr}(H^*_i H_i) = \text{tr} ( \sum_i  H_i^* H_i)  \leq d \|\sum_i H_i^* H_i \|.
    \end{equation}
   This, with $\theta<1$, implies that our expression is bounded by 
   $$
      \text{const} ~ d(x,y)^2  \sum_{\omega_1 , \dots , \omega_j} 
     \| R^{k-j}(1) (\psi_{\omega_j} \dots \psi_{\omega_1}x)\| 
     ~\| m^{(j-1)}(\psi_{\omega_{j-1}} \dots \psi_{\omega_1}y) \|^2$$
      $$   \leq  \text{const} ~ d(x,y)^2\sup_n\| R^{n}(1)\|_\infty \sum_{\omega_j}  \|R^{j-1}(1)(y)\|
     \leq  \text{const} ~ d(x,y)^2 (\sup_n\|R^n(1)\|)^2. $$
  \end{proof}
\end{lemma}

\begin{lemma} \label{key estimate 2}
There exist $D_1, D_2 >0$ such that 
$$
  \|p^*R^kh(x) - p^*R^kh(y) \|    
 \leq  D_1\|h\|_1 d( \psi_{\omega_k} \dots \psi_{\omega_1} x , 
 \psi_{\omega_k} \dots \psi_{\omega_1} y)  + 
\|h\|_\infty  D_2 d(x,y).
$$
for every $k \in \N$ and $\omega \in \Omega_k$
\begin{proof}
\begin{align*}
 & \|   p^*R^kh(x) -  p^*R^kh(y)  \|  \\ \leq  &    
\| \sum_{\omega \in \Omega_k} 
m^{(k)*} (\psi_{\omega_k} \dots \psi_{\omega_1} x) 
( h(p \psi_{\omega_k} \dots \psi_{\omega_1} x ) -   h(p \psi_{\omega_k} \dots \psi_{\omega_1} y ))
m^{(k)}( \psi_{\omega_k} \dots \psi_{\omega_1}x)\|\\
&   +
\|
 \sum_{\omega \in \Omega_k}  
m^{(k)*}(\psi_{\omega_k} \dots \psi_{\omega_1}x) 
h( p  \psi_{\omega_k} \dots \psi_{\omega_1} y ) 
(m^{(k)}( \psi_{\omega_k} \dots \psi_{\omega_1} x)- m^{(k)}( \psi_{\omega_k} \dots \psi_{\omega_1}y))\| \\
\\ & +
\|\sum_{\omega \in \Omega_k}  
(m^{(k)*}( \psi_{\omega_k} \dots \psi_{\omega_1}x)-m^{(k)*}( \psi_{\omega_k} \dots \psi_{\omega_1}y))  
 h(p \psi_{\omega_k} \dots \psi_{\omega_1} y ) 
m^{(k)}( \psi_{\omega_k} \dots \psi_{\omega_1}y)
\|.
\end{align*}
If we use Lemma \ref{CauchySchwarz} several times we get that this is bounded by 
\begin{align*}
&
\| \sum_{\omega \in \Omega_k} 
m^{(k)*}( \psi_{\omega_k} \dots \psi_{\omega_1} x) 
m^{(k)}( \psi_{\omega_k} \dots \psi_{\omega_1} x)\|\sup_{\omega \in \Omega_k}
 \|( h(p  \psi_{\omega_k} \dots \psi_{\omega_1} x ) -   h(p \psi_{\omega_k} \dots \psi_{\omega_1} y ))\|
\\
& +
\| \sum_{\omega \in \Omega_k}  
m^{(k)*}(  \psi_{\omega_k} \dots \psi_{\omega_1} x) m^{(k)*}( \psi_{\omega_k} \dots \psi_{\omega_1} x)  \|^{1/2} 
 \sup_{\omega \in \Omega_k}  \|h(p \psi_{\omega_k} \dots \psi_{\omega_1} y )\| \\
 &
\|  \sum_{\omega \in \Omega_k} (m^{(k)*}( \psi_{\omega_k} \dots \psi_{\omega_1} x)- m^{(k)*}( \psi_{\omega_k} \dots \psi_{\omega_1}  y))
 (m^{(k)}( \psi_{\omega_k} \dots \psi_{\omega_1} x)- m^{(k)}(  \psi_{\omega_k} \dots \psi_{\omega_1} y))\|^{1/2} \\
\\ + & 
\|\sum_{\omega \in \Omega_k}  
(m^{(k)*}(  \psi_{\omega_k} \dots \psi_{\omega_1}x)-m^{(k)*}( \psi_{\omega_k} \dots \psi_{\omega_1} y))
 (m^{(k)}( \psi_{\omega_k} \dots \psi_{\omega_1} x)-m^{(k)}( \psi_{\omega_k} \dots \psi_{\omega_1} y))  \|^{1/2} \\ 
&  \sup_{\omega \in \Omega_k}  \|
h( p \psi_{\omega_k} \dots \psi_{\omega_1} y )\|~
\|\sum_{\omega \in \Omega_k}   
m^{(k)*}( \psi_{\omega_k} \dots \psi_{\omega_1} y)m^{(k)}( \psi_{\omega_k} \dots \psi_{\omega_1}  y)
\|^{1/2}.
\end{align*}
This is bounded by 

  $$\|R^k 1 (py)\| ~ \sup_\omega
\|( h(p \psi_{\omega_k} \dots \psi_{\omega_1}  x ) -   h( p\psi_{\omega_k} \dots \psi_{\omega_1}y ))\|
+
 \| h \|_\infty  (\|R^k 1 (px)\|^{1/2} +\|R^k 1 (py)\|^{1/2}  ) $$
 $$
\|\sum_\omega 
(m^{(k)*}( \psi_{\omega_k} \dots \psi_{\omega_1} x)-m^{(k)*}(  \psi_{\omega_k} \dots \psi_{\omega_1} y))
(m^{(k)}( \psi_{\omega_k} \dots \psi_{\omega_1}  x)-m^{(k)}( \psi_{\omega_k} \dots \psi_{\omega_1}   y))  \|^{1/2}. 
$$

The claim follows from this and Lemma \ref{key estimate 1}.
\end{proof}
\end{lemma}

\begin{lemma}\label{equivalent norms}
The map
  $$
  h \mapsto \|h\|_\infty + \sup_{x\neq y} \frac{\|p^*h(x)- p^*h(y)\|}{d(x,y)}
  $$
  defines an equivalent norm on $L$.
  \begin{proof}
    A straightforward computation gives the following 
    estimate
  $$
     \sup_{ \| p^*s\|_1 \leq 1} \|hs\|_\infty + 
      \sup_{x \neq y} \frac{\| p^* h s (x) - p^* hs(y)\| }{d(x,y) }
      \leq  2 ( \| h\|_\infty  + \sup_{ x \neq y} 
      \frac{\| p^*h (x) - p^*h(y) \| }{d(x,y)}).
  $$
    On the other hand, we can find a finite cover of $X$ by open sets $V$ such that the restriction of $\xi$ to these sets is trivial. By the Lebsgue covering lemma there exists an $\epsilon > 0$ such  that 
     whenever $x,y\in\tilde X$ have $d(px,py) < \epsilon /2 $, there exists a $V$ in this cover such that $px,py\in V$. We can find then a unit vector  $v_{px,py} \in \C^d$
    such that 
    $$
       \| p^*h(x) - p^*h(y) \| = \|(p^*h(x) - p^*h(y))v_{px,py}\|.
    $$
    Moreover, for every $x,y\in\tilde X$ with $px,py \in V$, we can construct a section $s^{px,py} \in S_1$ such that 
     $p^*s_{px,py}(z) = v_{px,py}$ for every $z$ with $pz \in V$, and, by constructing a Lipschitz partition of unity subordinated to our cover, we can assume that $\|s_{px,py}\|_1$ is uniformly bounded by some constant $D$ that does not depend on $V,x,y$. 
     \par 
    Now 
    $$
          \frac{ \| p^*h(x) - p^*h(y) \|}{d(x,y)}  = 
          \frac{ \| (p^*h(x) - p^*h(y))v_{px,py}  \|}{d(x,y)} 
        \leq   \sup_{\| p^*s\|_1 \leq D}
       \frac{ \| p^*hs(x) - p^*hs(y)  \|}{d(x,y)} 
    $$$$\leq  \sup_{\| p^*s\|_1 \leq D} \sup_{x\neq y\in p\- V } 
       \frac{ \| p^*hs(x) - p^*hs(y)  \|}{d(x,y)}.
     $$
   If $d(px, py) \geq \epsilon/2$ then
   $$
    \frac{ \| p^*h(x) - p^*h(y) \|}{d(x,y)} \leq 4 \| h\|_\infty / \epsilon ,
   $$
    so we get a constant such that 
    $$
    \|h\|_\infty + \sup_{x\neq y} \frac{\|p^*h(x)- p^*h(y)\|}{d(x,y)} \leq D \| h \|_1.
    $$
  \end{proof}
\end{lemma}

\begin{theorem}
  The essential spectral radius of $R : L \rightarrow L$
  is less than  $\theta$. Since $\sup_k \| R^k h \|_\infty  < \infty$, the spectral radius is at most $1$.

  Moreover, in the main example with $r(px) = p(Ax)$
  where $A$ is  a linear map on $\R^n$, 
  the essential spectral radius is less than 
  the spectral radius of $A\-$.
  \begin{proof}
    By the Lemmas \ref{key estimate 2} and \ref{equivalent norms}
    we get $D_1, D_2 > 0$ such that 
    $$
    \| R^k h \|_1 \leq D_1 \theta^k \| h\|_1 + D_2 \| h \|_\infty ,
    $$
    in the general case, and 
    $$
    \| R^k h \|_1 \leq D_1 \| A^{-k}\| \| h\|_1 + D_2 \| h \|_\infty, 
    $$  
    in the case with the linear map $A$.

  By the Arzela-Ascoli theorem, bounded sets in $S_1$ are precompact in $S$. This implies  that 
  bounded sets in $L$ are precompact in $\text{End}(\xi)$.
The claim now follows from Corollary
    \ref{quasicompact} and the spectral radius formula. 
  \end{proof}
\end{theorem} 

\begin{theorem}\label{cesaro} The Cesaro means 
$k\- \sum_{j = 1}^k R^jg$ converges with respect to the $\| \cdot \|_\infty$-norm for every 
$g \in \text{End}(\xi)$. The map $T_1 : g \mapsto \lim_k k\- \sum_{j = 1}^k R^j   g$ 
  defines a completely positive idempotent acting on $\text{End}(\xi)$ such that $T_1|_L \in B(L)$ and 
$T_1 \text{End} (\xi ) = \{g \in \text{End}(\xi) | Rg = g\} = \{g \in L | Rg = g\}  $.
\end{theorem}
\begin{proof}
  Theorem \cite[VIII.5.1]{Dunford}, states that if 
  $R$ is a bounded operator on a Banach space $\mathcal{B}$ with uniformly bounded iterates, then 
  the set of $x \in\mathcal{B}$ such that the Cesaro means $n^{-1} \sum_{j= 1}^n R^j x$ converge is a closed subspace consisting of 
  all $y \in\mathcal B$ such that $n^{-1} \sum_{j= 1}^n R^jy$ is weakly sequentially compact and $\lim_n n^{-1} R^n y  = 0$.
  To apply this to our setting, note that the sequence $n^{-1} \sum_{j = 1}^n R^j g$,  $n \in \N$,   
  is  uniformly bounded with respect to the 
  $\| \cdot \|_1$-norm  whenever $g \in L$, so, using the Arzela-Ascoli theorem, it has a convergent 
  subsequence  with respect to the $\| \cdot \|_\infty $-norm. Moreover, $\lim_n  n^{-1} R^n g = 0$ 
  for every $g \in L$. Since $L \subset \text{End}(\xi)$ is dense with respect to the 
  $\| \cdot \|_\infty$-norm, this theorem implies that $n^{-1} \sum_{j= 1}^n R^jg$ converges,    
   $RT_1 g = T_1g$, and $Rg = g$ implies that $T_1g = g$ for every $g \in \text{End}(\xi)$.

  Since 
   the uniform limit of 
  a sequence of uniformly bounded Lipschitz continuous 
  sections  is also Lipschitz continuous  with Lipschitz constant less than 
  the bound, the limit of the convergent subsequence is in $L$ and the restriction of $T_1$ is bounded with respect to 
  the $\| \cdot \|_1$-norm. 

  Let $h \in \text{End}(\xi)$ such that $Rh = h$ and pick a 
  sequence $\{h_n\} \subset L$ such that 
  $\lim_n \| h_n - h \|_\infty = 0$. Now $T_1h_n \in T_1 L$ 
  converges to $h$ with respect to the $\| \cdot \|_\infty$-norm. 
  Moreover, we have that $T_1 |_L = \{g \in L | Rg = g\}$ is finite dimensional so it is closed 
  in both $L$ and $\text{End}(\xi)$. This implies that $h \in T_1 L$
  and $\{g \in \text{End} (\xi) | R g = g \} = \{ g \in L | Rg = g\}$. 
  
  A  map $P \in B(\text{End}(\xi))$ is completly positive if
  $ a \otimes h \mapsto a \otimes P(h)$ 
  defines a positive map $P^{(k)}$ on  $M_m(\C) \otimes \text{End}(\xi)$ 
  for every $k \in \N$. We see that 
  $h \mapsto m^*(p\psi \cdot    )   h(p \psi \cdot )    m( p \psi \cdot )$ 
  defines a  completely positive map. This implies that $R$ is also completely 
  positive. Since $T_1$ is the limit of Cesaro means of 
  the completely positive operators $R^k$, that converges in norm for every 
  $h \in \text{End}(\xi)$, we see that $T_1$ is also completely positive. Indeed, if $(h_{ij})_{i,j=1}^n$ is a positive matrix with coefficients in $\End\xi$, then $(R^kh_{ij})_{ij}$ is a positive matrix, and since each Cesaro mean $\frac{1}{k}\sum_{l=1}^k R^lh_{ij}$ converges to $T_1h_{ij}$ we get that $(T_1h_{ij})_{i,j}$ is a positive matrix, so $T_1$ is completely positive.
   
\end{proof}  

The essential spectral radius of similar transfer operators is analyzed in 
several other places in the literature. In \cite{BraJo02} this is computed with
the help of a theorem by Ionescu  and Tulcea \cite{MR0037469}, while in \cite{MR1016871} this is computed with techniques  from symbolic dynamics and dynamical zeta
functions.

\section{The  $C^*$-algebra of continuous harmonic maps} \label{c*al}
\begin{definition}
An element $h\in\End\xi$ is called {\it harmonic } with respect to the transfer operator $R$ if $Rh=h$. We denote by $\mathfrak H$ the set of all continuous harmonic functions $h\in\End\xi$. 
\end{definition}
Suppose $h$ is an invertible and 
strictly positive $R$-harmonic bundlemap.
We introduce the normalized 
transition operator $\tilde R$ on  the bundlemaps on $\xi$. Let 
$\tilde m = h^{1/2} m (h ^{-1/2}\circ r )$. Note that $\tilde m$ is Lipschitz because $h$ is Lipschitz (by Theorem \ref{cesaro}) and $\det h$ is bounded away from 0.
\par
Now $\tilde R$ is a completely  positive map 
that preserves the identity with a corresponding  
 completely positive idempotent  $\tilde T_1$.   A theorem due to Choi and Effros \cite{ChoiEffr} states that 
 the image  of a 
completely positive map $T$ that preserves the identity, 
equipped with the usual norm,  $*$-operation and the product 
$a , b \mapsto  T(ab)$ form a $C^*$-algebra, i.e. the set of 
$\tilde R$-harmonic maps is  a $C^*$-algebra when  equipped with usual 
$*$-operation, the usual norm and the product $a, b \mapsto \tilde T_1 (ab)$.
  
The map $g \mapsto h^{-1/2} g h^{-1/2}$ is a 
linear bijection from the $ R$-harmonic bundle maps to the 
$\tilde R$-harmonic bundlemaps. This means that 
the $R$-harmonic bundlemaps with the usual $*$-operation, 
the norm $\| a \|_h = \| h^{-1/2} a h ^{-1/2} \|$ and the 
product $a, b \mapsto h^{1/2}\tilde T_1 ( h^{-1/2} a h\- b  h^{-1/2})h^{1/2} = 
T_1( a h \- b)$ form a $C^*$-algebra. Let $\mathfrak H$ denote this algebra and define $a*b = T_1(a h\- b)$.

Since the essential spectral radius of $R$ is strictly less that $1$, 
$\mathfrak H$ is finite dimensional, so we get the following:
\begin{theorem}
 The set $\mathfrak H$ of continuous harmonic functions for the transfer operator $R$ forms a finite dimensional $C^*$ algebra with the usual addition, adjoint and multiplication by scalars and with the product defined by $h_1*h_2:=T_1(h_1h\-h_2)$ and norm $\|f\|_h:=\|h^{-1/2}fh^{-1/2}\|$.
 \end{theorem}
Every finite dimensional $C^*$-algebra is isomorphic 
to $M_{k_1} (\C) \oplus \dots \oplus M_{k_r}(\C)$ for $k_1, \dots k_r \in \N$. See for instance 
Theorem  \cite[III.1.1]{Davidson} or \cite[ 7.1]{Rordam}.
In what follows, we will give a partial description of this algebra in some situations.

\section{Refinement operators}\label{refi}
 We define here the refinement operator $M$. This operates at the level of $\tilde X$, the simply connected covering space of $X$. On $\tilde X$ we have a unitary ``dilation operator'' $U$. It is defined by the following:
\begin{proposition}
 The measure $\tilde\mu$ satisfies the following quasi-invariance equation:
$$
\frac{ d \tilde \mu \circ \tilde r \-}{d \tilde \mu } = q.
$$
 The operator $U$ defined on $\H:=L^2(\tilde X,\tilde\mu)$ by $Uf=q^{1/2}f\circ\tilde r$ is unitary, with $U^*f=q^{-1/2}f\circ\tilde r^{-1}$. \par
 Define $\pi(a)f (x)  = a \circ p (x) f(x)$ for $a\in C(X)$ and $f\in\H$. Then $\pi$ is a 
representation of $C(X)$ in $B(\H)$, and $( \H, \pi, U)$ 
form a covariant representation of the dynamical system $r : X \rightarrow X$ in the sense
that $U \pi(f) U^* = \pi(f \circ r)$ for every $f \in C(X)$.
\end{proposition}

\begin{proof}
Using the stong invariance \eqref{eqstrinv} of the measure $\mu$, equation \eqref{eqpermu} and \eqref{eqpsi}, we have for all $f\in C_c(\tilde X)$:
  \begin{align*}
    \int_{\tilde X} f d \tilde \mu & = \int_{\tilde X} f \circ \tilde r 
    \frac{ d \tilde \mu \circ \tilde r \-}{d \tilde \mu } d \tilde \mu  
    = \int_X \sum_{g \in G} (f \circ \tilde r   
    \frac{ d \tilde \mu \circ \tilde r \-}{d \tilde \mu } ) \circ g d \mu 
     =  \int_X q\-  \sum_{i=1}^q \sum_{g \in G} (f \circ \tilde r   
    \frac{ d \tilde \mu \circ \tilde r \-}{d \tilde \mu } ) \circ g \circ \psi_i d \mu \\
    & =\int_Xq\-\sum_{i=1}^q\sum_{g\in G}(f\circ\tilde r\frac{ d \tilde \mu \circ \tilde r \-}{d \tilde \mu } )\tilde r^{-1}(Ag)g_i\,d\mu= \int_X q\-  \ \sum_{g \in G} (f    
    \frac{ d \tilde \mu \circ \tilde r \-}{d \tilde \mu }\circ \tilde r \-  ) \circ g  d \mu  
     \\ & = \int_{\tilde X} q\- f    
    \frac{ d \tilde \mu \circ \tilde r \-}{d \tilde \mu }\circ \tilde r \-    d \tilde \mu.
  \end{align*}
    \par
  The other statements follow from some easy computations.
\end{proof}

\begin{definition}
As in \cite{PackerRieffel1}, we define 
$\Xi
$
 to be $f \in C_b(\tilde X)$ (bounded continuous functions on $\tilde X$) such that $\sum_{g \in G} |f(gx)|^2$ is bounded and continuous 
 for every $x \in \tilde X$. Following their proof, we see that when 
 $\Xi \subset C_b(\tilde X)$ is equipped with the inner product
 $$\langle f_1, f_2 \rangle' =  \sum_{g \in G} \cj f_1 \circ g    f_2 \circ g,$$  
 we have a  $C(X)$-Hilbert module.
\par
Let $\text{Hom}_{C(X)}(S,\Xi)$ be the set of adjointable $C(X)$-linear maps between the Hilbert modules $S$ and $\Xi$. 
 
\par
The {\it refinement operator} associated to the map $m$,
$$M : \text{Hom}_{C(X)}(S, \Xi) \rightarrow \text{Hom}_{C(X)}(S, \Xi)$$ is defined by
$$M W s = U\- W (q^{-1/2}m s \circ r)=(Wms\circ r)\circ\tilde r^{-1},\quad(W\in\text{Hom}_{C(X)}(S,\Xi),s\in S).$$
\end{definition}

\begin{proposition}
  If $W \in \text{Hom}_{C(X)}(S, \Xi)$ then $U W S \subset WS$ if and only if 
  there exists a bundlemap  $\tilde m : r^*\xi \rightarrow \xi$ such that $W$ is 
  a fixed point for the refinement operator $\tilde M$ associated to $\tilde m$. 
  \end{proposition}
  \begin{proof}
    If $ W \in \text{Hom}_{C(X)}(S, \Xi)$ such that $\tilde M W = W$, we see imediately from the
    definition of $M$ that $U W S \subset WS$. Conversely, suppose  $U W S \subset W S$. 
    Let $\zeta$ be  another 
    $C(X)$-module such that $\xi \oplus \zeta \simeq X \times \C^k $ (see \cite[Corollary 1.4.14]{Atiyah1}). Let $s_1, \dots , s_k$ 
    be linearly independent sections in $\xi \oplus \zeta$. The sections
    $\text{Pr}_1 s_j$,  $1 \leq j \leq k  $ generate $S$ as a $C(X)$-module. 
    There 
    exist $c_{i,j} \in C(X)$,
    $1 \leq i,j \leq k$ such that 
    $ U W \text{Pr}_1 s_j  = \sum_i  c_{i,j} \text{Pr}_1 W s_i $ for every $j$.
        Moreover $r^*\xi \oplus r^*\zeta \simeq r^*( \xi \oplus \zeta)$ 
    and 
    $$
    x \mapsto (x, \text{Pr}_1 s_i(rx)) \oplus (x, \text{Pr}_2 s_i(rx)),\quad(1 \leq i \leq k),
    $$
    form $k$ linearly independent sections 
    $s_1', \dots , s_k'$ in the bundle in $r^*\xi \oplus r^* \zeta$.
    Now define  $C : r^* \xi \oplus r^* \zeta \rightarrow \xi \oplus \zeta$
    such that $Cs_i' = \sum_j c_{i,j} s_j$,
    and 
    $\tilde m = q^{1/2}\text{Pr}_1 C \text{Pr}_1$. This gives us that 
    $U\- W (q^{-1/2}\tilde m s \circ r) = W s$ for every $s \in S$.    
  \end{proof}

\begin{proposition}
  If  $W \in \text{Hom}(S, \Xi)$ such that
  \begin{enumerate}
    \item $M W = W$;
    \item $W^*W$ is an idempotent (Note that here $W^*$ is the adjoint of $W$ as a Hilbert module map between $S$ and $\Xi$ !);
    \item There exists an $f \in WS$ such that $f(\tilde x_0) \neq 0$;
  \end{enumerate}
  then $\{U^k WS \}_{k \in \Z}$ form a projective multiresolution analysis, i.e. 
  the following conditions are satisfied
  \begin{enumerate}
  \item $WS$ is a projective  submodule in $\Xi$;
  \item $WS \subset U\-WS$;
  \item $\cup_{k \in \Z} U^k WS$ is dense in $\Xi$;
  \item $\cap_{k \in \Z} U^k WS = \{0\}$;
  \end{enumerate}
  \begin{proof}
    If $W^*W$ is a projection, then $WW^*$ is also a projection, $WS = W^*W \Xi$
    and $W |_{W^*WS}$ is an isometry onto its image, \cite{Lance}. This implies that $WS$ is isomorphic to 
    the image of an idempotent on a projective module and therefore it must be projective also. 
    By proposition 13 and 14 in \cite{PackerRieffel1}, we get conditions 2, 3 and 4. 
  \end{proof}
\end{proposition}

We 
have the following important intertwining relation between the refinement operator and the transfer operator:
\begin{theorem}\label{thinter}
If $ W_1, W_2  \in  \text{Hom}_{C(X)}(S, \Xi)$, then 
  $$
  R( W_1^* W_2) = (MW_1)^* MW_2 .
  $$
  In particular, if $MW_1=W_1$ and $MW_2=W_2$ then $R(W_1^*W_2)=W_1^*W_2$, i.e., $W_1^*W_2$ is harmonic.
  \begin{proof}
   Let $s_1, s_2 \in S$
    \begin{align*}
       \langle M W_1 s_1, M W_2 s_2 \rangle' \circ p =&\sum_{g\in G}\cj {(W_1ms_1\circ r)\circ\tilde r^{-1}\circ g } 
       (W_2ms_2\circ r)\circ\tilde r^{-1} \circ g \\  
       = &   \sum_{i=1}^q \sum_{g \in G} \cj{ ( W_1 m s_1 \circ r ) \circ \tilde r \- \circ Ag \circ g_i}
      ( W_2 m s_2 \circ r ) \circ \tilde r \- \circ Ag \circ g_i    \\
       = &   \sum_{i=1}^q \sum_{g \in G} \cj{ ( W_1 m s_1 \circ r )  \circ g \circ \psi_i } 
      ( W_2 m s_2 \circ r ) \circ g \circ \psi_i    \\
        = &   \sum_{i=1}^q  \langle  W_1 m s_1 \circ r , W_2 m s_2 \circ r \rangle' \circ p \circ \psi_i \\ 
	= &  \sum_{i=1}^q \langle s_1 \circ r , (m^* W_1^* W_2 m) s_2 \circ r \rangle \circ p \circ \psi_i \\
	= & \langle s_1, R(W_1^*W_2) s_2 \rangle \circ p. 
    \end{align*}
  \end{proof}
\end{theorem}

This suggests a method to construct embeddings $W : S \rightarrow \Xi$,  
such that $U W S \subset W S$. Pick an $m : r^* \xi \rightarrow \xi$ and 
let $R$ be the corresponding transfer operator with a basis of  
$R$-invariant functionals on $\text{End}(\xi)$,  $\tau_1, \dots, \tau_r$,
and a positive and invertible harmonic map $h$. 
Then look for a fixed point $W$ of the refinement operator $M$,  such that 
$\tau_j(h ) =  \tau_j(W^*W)$ for every $1 \leq j \leq r$. Then $W^*W = h$ and 
$W$ is injective. Moreover, if $h = 1$ we also know that $W$ is an isometry. 

This method is essentially the same as a technique often encountered in the wavelet literature. 
With some assumptions on the low-pass filter, the infinte product expansion 
yields a scaling function  that generates a multiresolution analysis
in $L^2(\R)$. This scaling function has orthonormal $\Z$-translates if and only if the constants are the 
only fixed points for the corresponding transfer operator, see for instance \cite{Law90}.  

Suppose $W : S  \rightarrow  \Xi$ is of the form  
$W s = \langle s_0, s \rangle f$, with $s_0\in S $ and $f\in  \Xi$, then $W$ is adjointable, and a direct computation gives the following two identities:
    \begin{equation}\label{eqmkw}
        M^k  W s (x)  =  \langle  s_0 \circ p( \tilde r^{-k} x)    
      ,  m^{(k)}  ( \tilde r^{-k} x) s \circ p (x) \rangle 
      f (  \tilde r^{-k}x).
    \end{equation}
    
    \begin{equation}\label{eqmk*w}    
      (M^kW)^* f \circ p(x)  
      =   \sum_{\omega \in \Omega_k}  
      m^{(k)*}( \psi_{\omega_k} \dots \psi_{\omega_1}x) (W^*f\circ \tilde r^k) \circ 
      p( \psi_{\omega_k} \dots \psi_{\omega_1}x).
    \end{equation}
    Moreover, there exist $M_1, M_2 > 0$ such that  
    $\|M^kW \| \leq M_1 \|W \|$ and $\|(M^k W)^*\| \leq M_2 \| W^* \|$ for every $k \geq 0$. 
\begin{definition}
Let $\Xi_1$ be the set of $f \in \Xi$ such that there exists a $D > 0$ such that
$$
(\sum_{g \in G} | f( g x) - f(g y) |^2)^{1/2} \leq D d(x,y),
$$
for every $x,y \in \tilde X$.
\end{definition}
\begin{lemma}
  $\Xi_1$ is a $Lip_1(X)$-module and $f_1, f_2 \in \Xi_1$ implies that 
$\langle f_1, f_2 \rangle' \in \text{Lip}_1(X)$.
\end{lemma}
  \begin{proof}
    $$
       | \langle f_1, f_2 \rangle '(x) - \langle f_1, f_2 \rangle'(y) | 
      \leq  |\sum_{g \in G} \cj f_1 (gx) ( f_1(gx) - f_2(gy) ) | 
      +  |\sum_{g \in G}   \cj{(f_1(gx) - f_2(gy))} f_2(gy)   |
      $$$$
       \leq  ( \langle f_1, f_1 \rangle'(x))^{1/2} (\sum_{g \in G} | f_2( g x) - f_2(g y) |^2)^{1/2}
+ ( \langle f_2, f_2 \rangle'(y))^{1/2} (\sum_{g \in G} | f_1( g x) - f_1(g y) |^2)^{1/2}.
    $$
    
  \end{proof}

\begin{lemma}\label{lemlipmk}
  Let $W \in \text{Hom}_{C(X)}(S, \Xi)$ such that
$W s (x) = \langle s'(px) , s(px)\rangle f'(x)$  with $s'\in S_1$ and
$f'\in\Xi_1$.
   There exist  $D_1, D_2  > 0$
  such that 
  $$
  ( \sum_{g \in G} | (M^k W s)(gx ) - (M^kWs)(gy) |^2)^{1/2} \leq D_1 \|s\|_\infty d(x,y) +  D_2 \| s(px) - s(py)\|,
  $$
  for every $x,y \in\tilde X$.
  \end{lemma}
  \begin{proof}
We have
    \begin{align}
        &( \sum_{g \in G} | (M^k W s)(gx ) - (M^kWs)(gy) |^2)^{1/2} \\
      \leq & (\sum_{g \in G}  |   \langle  s' \circ p( \tilde r^{-k} gx)-  s' \circ p( \tilde r^{-k} gy)    
      ,  m^{(k)}  ( \tilde r^{-k} gy) s \circ p (y) \rangle 
      f' (  \tilde r^{-k}gy) |^2)^{1/2}  \label{eq1}\\
      & +  (\sum_{g \in G}  |   \langle  s' \circ p( \tilde r^{-k} gx)    
      , (  m^{(k)}  ( \tilde r^{-k} gx) -  m^{(k)}  ( \tilde r^{-k} gy)   )s \circ p (y) \rangle 
      f' (  \tilde r^{-k}gy) |^2)^{1/2}\label{eq2} \\
         & +  (\sum_{g \in G}  |   \langle  s' \circ p( \tilde r^{-k} gx)    
      ,   m^{(k)}  ( \tilde r^{-k} gx) ( s \circ p(x) - s \circ p (y) \rangle 
      f' (  \tilde r^{-k}gy) |^2)^{1/2}\label{eq3} \\
   & +  (\sum_{g \in G}  |   \langle  s' \circ p( \tilde r^{-k} gx)    
      ,   m^{(k)}  ( \tilde r^{-k} gx) s \circ p (x) \rangle 
      (f' (  \tilde r^{-k}g x) - f' (  \tilde r^{-k}gy)) |^2 )^{1/2}. \label{eq4}
    \end{align}
    We estimate the last four summands separately. \eqref{eq1} is bounded by 
    \begin{align*}
      (\sum_{g \in G} & \sum_{\omega \in \Omega_k} 
      \|  \  s' \circ p(  \psi_{\omega_k} \dots  \psi_{\omega_1} x)-  
      s' \circ p( \psi_{\omega_k} \dots \psi_{\omega_1} y) \|^2\cdot \\        
      & \cdot \|  m^{(k)}  (  \psi_{\omega_k} \dots \psi_{\omega_1} y )s \circ p (y)\|^2 ~ | 
      f' (  g \psi_{\omega_k} \dots \psi_{\omega_1} y) |^2)^{1/2} \\
       =  & (\sum_{\omega \in \Omega_k} 
      \|  \  s' \circ p(  \psi_{\omega_k} \dots  \psi_{\omega_1} x)-  
      s' \circ p( \psi_{\omega_k} \dots \psi_{\omega_1} y) \|^2 \\ &  
     \|  m^{(k)}  (  \psi_{\omega_k} \dots \psi_{\omega_1} y )s \circ p (y)\|^2 ~ \langle 
      f', f'   \rangle' ( \psi_{\omega_k} \dots \psi_{\omega_1} y )  )^{1/2} \\
      & \leq \text{const}\,  \theta^k d(x,y) d \| R^k 1 \|_\infty \|s \|_\infty \| \langle f' , f' \rangle' \|_\infty^{1/2}. 
      \end{align*}
      (We used \eqref{eqhi} for the last estimate.)\par
\eqref{eq2} is bounded by 
    \begin{align*}
     (\sum_{g \in G}  \sum_{\omega \in \Omega_k} &
      \|  \  s' \circ p(  \psi_{\omega_k} \dots  \psi_{\omega_1} x) \|^2\cdot \\ 
          \|  &( m^{(k)}  (  \psi_{\omega_k} \dots \psi_{\omega_1} y )    
       -  m^{(k)}  (  \psi_{\omega_k} \dots \psi_{\omega_1} x ) ) s \circ p (y)\|^2 ~ 
      | f' (  g \psi_{\omega_k} \dots \psi_{\omega_1} y) |^2
           )^{1/2} \\
      \leq & \| s \|_\infty  \|s'\|_\infty  \| \langle f', f' \rangle' \|^{1/2}_\infty      
       (\sum_{\omega \in \Omega_k} \| m^{(k) } ( \psi_{\omega_k } \dots \psi_{\omega_1} x)
      -  m^{(k) } ( \psi_{\omega_k } \dots \psi_{\omega_1} y) \|^2)^{1/2} \\
      \leq & \text{const}\, \|s\|_\infty ~d(x,y).
      \end{align*}
    and we used Lemma \ref{key estimate 1} in the last inequality. \par
    Moreover, \eqref{eq3}
    is bounded by 
    \begin{align*}
      (\sum_{g \in G}  \sum_{\omega \in \Omega_k}& 
      \|  \  s' \circ p(  \psi_{\omega_k} \dots  \psi_{\omega_1} x) \|^2 
        \|   m^{(k)}  (  \psi_{\omega_k} \dots \psi_{\omega_1} y )    
       (  s \circ p (x) - s \circ p(y) ) \|^2 ~ 
      | f' (  g \psi_{\omega_k} \dots \psi_{\omega_1} y) |^2
           )^{1/2} \\
       \leq & \| s' \|_\infty d \| R^k 1 \|_\infty \| s(px) - s(py )\|  \| \langle f' , f' \rangle' \|^{1/2}. 
    \end{align*}
    Finally, \eqref{eq4} is bounded by the following expression
    $$
        (\sum_{g \in G}  \sum_{\omega \in \Omega_k} 
      \|  \  s' \circ p(  \psi_{\omega_k} \dots  \psi_{\omega_1} x) \|^2   
         \|   m^{(k)}  (  \psi_{\omega_k} \dots \psi_{\omega_1} y )    
         s \circ p (x) \|^2 ~ 
      |f' (  g \psi_{\omega_k} \dots \psi_{\omega_1} x) -  f' (  g \psi_{\omega_k} \dots \psi_{\omega_1} y) |^2
           )^{1/2} $$
           $$
      \leq  \| s' \|_\infty \|s \|_\infty  ( \sum_{\omega \in \Omega_k} 
      \| m^{(k)}  ( \psi_{\omega_k} \dots \psi_{\omega_1}x    ) \|^2  )^{1/2}     
      \,\text{const} \,
      d(\psi_{\omega_k} \dots \psi_{\omega_1} x
	 ,\psi_{\omega_k} \dots \psi_{\omega_1} y) 
     $$
     $$ \leq \,   \text{const}\,\|s\|_\infty\, \| R^k 1\|_\infty \theta^k d(x,y).
    $$
 \end{proof}

\section{An example: low-pass filters}\label{exam}
In this section we will assume that $m$ satisfies a certain low-pass condition. This will imply that the iterates of the refinement operator will converge to some fixed points, which correspond to scaling functions. These in turn will generate harmonic functions which are minimal projections in the algebra $\mathfrak H$.
\par 
To simplify the notation we will identify $m\circ p$ and $p^*m$ with $m$.
\begin{definition}\label{defel}
Let $1\leq l\leq d$. We say that a $d\times d$ matrix $a$ satisfies the $E(l)$ condition if $\|a\|=1$, $1$ is the only eigenvalue of $a$ on the unit circle and both the geometric and algebraic multiplicities of $1$ are equal to $l$.
\end{definition}

In the next theorem, we show how certain ``scaling functions'' $\mathcal W_v$ can be constructed from each eigenvector $v$ of $m(x_0)$ corresponding to the eigenvalue $1$. 

\begin{theorem}\label{prop5-1}
Assume $R1=1$ and $m(x_0)$ satisfies the $E(l)$ condition. 
\begin{enumerate}
\item The limit 
$$\mathcal{P}(\tilde x):=\lim_{k\rightarrow\infty}m^{(k)}(r^{-k}(\tilde x)),\quad (\tilde x\in\tilde X)$$
exists for all $\tilde x\in\tilde X$, and is uniform on compact sets. 
The following refinement equation is satisfied:
\begin{equation}\label{eqscalem}
\mathcal P(\tilde r^{-1}\tilde x)m(\tilde r^{-1}\tilde x)=\mathcal P(\tilde x),\quad(\tilde x\in\tilde X).
\end{equation}
$\mathcal P(\tilde x_0)$ is the orthogonal projection onto the eigenspace $E_1$ of $ m(x_0)$ that corresponds to the eigenvalue $1$. The range of $\mathcal P(\tilde x)$ is contained in $E_1$, i.e.
\begin{equation}\label{eqpinv}
 m(x_0)\mathcal P(\tilde x)=\mathcal P(\tilde x),\quad(\tilde x\in\tilde X).
\end{equation}
For all $g\in G$, $g\neq 1$, and all $v\in E_1$,
\begin{equation}\label{eqpgx-0}
\mathcal P(g\tilde x_0)v=0.
\end{equation}
\item Fix a function $f\in\Xi_1$ with $\ip{f}{f}'=1$ and $f(g\tilde x_0)=0$ for all $g\in G$, $g\neq 1$.
Let $s_0\in S_1$ with $s_0(x_0)=:v$. Define $W\in Hom(S,\Xi)$, $W_{s_0}s=\ip{s_0}{s}f$. Then, for all $s\in S$, $M^kW_{s_0}s$ converges uniformly on compact sets to a continuous $\mathcal{W}_vs\in\Xi$, $\mathcal W_vs(\tilde x)=\ip{v}{\mathcal{P}(\tilde x)s(p\tilde x)}$, and this defines a $\mathcal{W}_v\in Hom(S,\Xi)$. Moreover $\mathcal{W}_v$ depends only on the value $v$ of $s_0$ at $x_0$, not on the entire section $s_0$. In addition $\mathcal W_v$ satisfies the refinement equation 
\begin{equation}\label{eqscal-w}
M\mathcal W_v=\mathcal W_v
\end{equation}
\end{enumerate}

Also, $\mathcal W_vs(\tilde x_0)=\ip{\mathcal{P}(\tilde x_0)v}{s(x_0)}f(\tilde x_0)$.
\end{theorem}

\begin{proof}
(i)
First note that, since the vector bundle $p^*\xi$ is trivial, the terms of the product defining $m^{(k)}$ act on the same vector space so it makes sense to talk about the convergence of this product (see also Remark \ref{remgraph}).
\par
Since $ m(x_0)$ satisfies the $E(l)$ condition, we can find an invertible matrix $u$ such that $J:=u^{-1} m(x_0)u$ is in Jordan canonical form, with the $l\times l$ leading principal submatrix of $u^{-1} m(x_0)u$ being the identity $l\times l$ matrix. The other Jordan blocks correspond to eigenvalues $\lambda$ with $|\lambda|<1$.
\par
Then one can see that $J^p$ converges to the matrix $\left[\begin{array}{cc}I_l&0\\ 0&0\end{array}\right]$. This shows in particular that $\|m(x_0)^p-m(x_0)^{p+p'}\|$ can be made arbitrarily small.
\par
Take now a compact subset $\tilde K$ of $\tilde X$. Since $\tilde r^{-1}$ is contractive towards the fixed point $\tilde x_0$, there are some constants $C>0$, $0<\theta<1$ depending only on the set $\tilde K$ such that $d(\tilde r^{-k}\tilde x,\tilde x_0)\leq C\theta^{k}$ for all $\tilde x\in\tilde K$ and $k\geq 0$.
\par
We have for $k\geq 0$, $p,p'\geq 0$
$$\|m^{(k+p)}(\tilde r^{-(k+p)}\tilde x)-m^{(k+p+p')}(\tilde r^{-(k+p+p')}\tilde x)\|=$$$$
\|m^{(p)}(\tilde r^{-(k+p)}\tilde x)m^{(k)}(\tilde r^{-k}\tilde x)-
m^{(p+p')}(\tilde r^{-(k+p+p')}\tilde x)m^{(k)}(\tilde r^{-k}\tilde x)\|$$
$$\leq 
\|m^{(k)}(\tilde r^{-k}\tilde x)\|\left\|m^{(p)}(\tilde r^{-(k+p)}\tilde x)-m^{(p+p')}(\tilde r^{-(k+p+p')}\tilde x)\right\|$$

\par
Since $R^k1=1$ we have that $\|m^{(k)}(\tilde r^{-(k)}\tilde x)\|\leq 1$. 
On the other hand, we have with Lemma \ref{key estimate 1}
$$\left\|m^{(p)}(\tilde r^{-(k+p)}\tilde x)-m^{(p+p')}(\tilde r^{-(k+p+p')}\tilde x)\right\|
\leq\|m^{(p)}(\tilde r^{-(k+p)}\tilde x)-m^{(p)}(\tilde r^{-p}\tilde x_0)\|
+\|m(x_0)^p-m(x_0)^{p+p'}\|$$$$
+\|m^{(p+p')}(\tilde r^{-(p+p')}\tilde x_0)-m^{(p+p')}(\tilde r^{-(k+p+p')}\tilde x)\|$$
$$\leq Dd(\tilde r^{-k}\tilde x,\tilde x_0)+\|m(x_0)^p-m(x_0)^{p+p'}\|+Dd(\tilde x_0,\tilde r^{-k}\tilde x)
\leq 2DC\theta^k+\|m(x_0)^p-m(x_0)^{p+p'}\|.$$
This proves that the sequence is uniformly Cauchy, therefore it is uniformly convergent.
\par
The scaling equation (\ref{eqscalem}) follows from the convergence of the product. 
Now let us compute $\mathcal P(\tilde x_0)$. We have $\mathcal P(\tilde x_0)=\lim_k m(x_0)^k$. Then  $\mathcal P(\tilde x_0)=\lim_km(x_0)^{2k}=\lim_km(x_0)^km(x_0)^k=\mathcal P(\tilde x_0)\mathcal P(\tilde x_0)$. 
\par
If $v\in\bc^d$, $ m(x_0)\mathcal P(\tilde x_0)v=\lim_km(x_0)^{k+1}v=\mathcal P(\tilde x_0)v$, so the range of $\mathcal P(\tilde x_0)$ is contained in the eigenspace $E_1$ of $ m(x_0)$ corresponding to the eigenvalue $1$. 
We have the following lemma:
\begin{lemma}\label{lemeigmx-0}
The eigenspaces corresponding to the eigenvalue $1$ for $ m(x_0)$ and $ m^*(x_0)$ are the same.
\end{lemma}
\begin{proof}
Let $v\in\bc^n$ such that $ m(x_0)v=v$. We may assume $\|v\|=1$. Since $R1=1$, we have $m^*(x_0)m(x_0)\leq 1$ so $\| m^*(x_0)\|\leq 1$. Using the Schwarz inequality we obtain
$$1=\ip{ m(x_0)v}{v}=\ip{v}{ m^*(x_0)v}\leq \|v\|\| m^*(x_0)v\|\leq 1,$$
therefore we have equalities in all these inequalities. This implies that $ m^*(x_0)v=\lambda v$ for some $\lambda\in\bc^n$, and $\lambda$ has to be $1$. Hence $ m^*(x_0)v=v$. The reverse implication can be proved analogously.  
\end{proof}
With Lemma \ref{lemeigmx-0} we have that $\mathcal P(x_0)^*=\lim_km^*(x_0)^k$ leaves $E_1$ fixed, so $\mathcal P(\tilde x_0)^*\mathcal P(\tilde x_0)=\mathcal P(\tilde x_0)$. But this implies that $\mathcal P(\tilde x_0)$ is a positive operator, and since it is also an idempotent, it must be equal to the orthogonal projection onto its range $E_1$.
\par
To prove equation (\ref{eqpinv}) we estimate, for $\tilde x\in\tilde X$:
$$\| m(x_0)\mathcal P(\tilde x)-\mathcal P(\tilde x)\|\leq \| m(x_0)\mathcal P(\tilde x) - m(x_0)m(p\tilde r^{-k}\tilde x)\dots m(p\tilde r^{-1}\tilde x)\|
+$$$$\| m(x_0)m(p\tilde r^{-k}\tilde x)\dots m(p\tilde r^{-1}\tilde x)-m(p\tilde r^{-(k+1)}\tilde x)\dots m(p\tilde r^{-1}\tilde x)\|+
\|m(p\tilde r^{-(k+1)}\tilde x)\dots m(p\tilde r^{-1}\tilde x)-\mathcal P(\tilde x)\|$$
But 
$$\| m(x_0)m(p\tilde r^{-k}\tilde x)\dots m(p\tilde r^{-1}\tilde x)-m(p\tilde r^{-(k+1)}\tilde x)\dots m(p\tilde r^{-1}\tilde x)\|\leq$$$$ \| m(x_0)- m(p\tilde r^{-(k+1)}\tilde x)\|\|m^{(k)}(\tilde r^{-k}\tilde x)\|\leq \| m(x_0)- m(p\tilde r^{-(k+1)}\tilde x)\|,$$
since $R1=1$ and $\| m(\tilde x)\|\leq 1$ for all $\tilde x$.
\par
Therefore all the terms can be made as small as we want, provided $k$ is big, and this proves (\ref{eqpinv}).
\begin{lemma}\label{lemmgx-0}
For all $v$ in the eigenspace $E_1$, $$ m(\psi_k\tilde x_0)v=0,\quad(2\leq
k\leq q).$$
\end{lemma}

\begin{proof}
Since $R1=1$ one has
$$\ip{v}{v}=\sum_{k=1}^q\ip{v}{m^*(p\psi_k\tilde x_0)m(p\psi_k\tilde
x_0)v}=\ip{m(x_0)v}{m(x_0)v}+\sum_{k=2}^q\ip{m(p\psi_k\tilde
x_0)v}{m(p\psi_k\tilde x_0)v}=$$$$\ip{v}{v}+\sum_{k=2}^q\|m(p\psi_k\tilde
x_0)v\|^2.$$
This implies the lemma.
\end{proof}
\par
Take now $g\in G$ , $g\neq 1$. Let $k$ be the first positive integer such
that $\tilde r^{-k}g\tilde x_0\not\in G\tilde x_0$. Then $p(\tilde
r^{-l}g\tilde x_0)=x_0$ for $1\leq l\leq k-1$ and $r(p(\tilde
r^{-k}g\tilde x_0))=x_0$, but
$p(\tilde r^{-k}g\tilde x_0)\neq x_0$. Using Lemma \ref{lemmgx-0} we
obtain for $j\geq k$
$$m^{(j)}(\tilde r^{-j}g\tilde x_0)v=m(p\tilde r^{-j}g\tilde
x_0)...m(p\tilde r^{-k}g\tilde x_0) m(x_0)\cdots m(x_0)v=0.$$
So $\mathcal P(g\tilde x_0)v=0$.
\par
(ii) With \eqref{eqmkw} we have
$$M^kW_{s_0}s(p\tilde x)=\ip{s_0\circ p(\tilde r^{-k}\tilde
x)}{m^{(k)}(\tilde r^{-k}\tilde x)s\circ p(\tilde x)}f(\tilde r^{-k}x)$$
Each term in this formula is uniformly convergent on compact sets so $M^kW_{s_0}s$ is uniformly convergent on compact sets to $\ip{s_0(x_0)}{\mathcal{P}(\tilde x)s\circ p(\tilde x)}$. Note that $f(\tilde x_0)=1$ because $\ip{f,f}'=1$ and 
$f(g\tilde x_0)=0$ for $g\neq 1$.
Using Fatou's lemma, and the intertwining relation in Theorem \ref{thinter}, we have
$$\ip{\mathcal W_vs}{\mathcal W_vs}'(p\tilde x)=\sum_{g\in G}\ip{\mathcal W_vs}{\mathcal W_vs}(g\tilde x)\leq\liminf_k\sum_{g\in G}\ip{M^k W_{s_0}s}{M^kW_{s_0}s}(g\tilde x)=$$
$$\liminf_k\ip{M^kW_{s_0}s}{M^kW_{s_0}s}'(p\tilde x)=
\liminf_k\ip{s}{(M^kW_{s_0})^*(M^kW_{s_0})s}(p\tilde x)=$$$$\liminf_k\ip{s}{R^k(W_{s_0}^*W_{s_0})s}(p\tilde x)\leq\ip{s}{s}\sup_{x\in X}\|W_{s_0}^*W_{s_0}(x)\|,$$
and we used $R^k(W_{s_0}^*W_{s_0})\leq R^k(\sup_{x}\|W_{s_0}^*W_{s_0}\|)=\sup_x\|W_{s_0}^*W_{s_0}\|$ in the last inequality ($R1=1$).
\par
Moreover, we have, with Lemma \ref{lemlipmk}
$$|\ip{\mathcal W_vs}{\mathcal W_vs}'^{\frac12}(p\tilde x)-\ip{\mathcal W_vs}{\mathcal W_vs}'^{\frac12}(p\tilde y)|
= |(\sum_{g\in G}\|\mathcal W_vs(g\tilde x)\|^2)^{\frac12}-(\sum_{g\in G}\|\mathcal W_vs(g\tilde y)\|^2)^{\frac12}|$$
$$\leq(\sum_{g\in G}\|\mathcal W_vs(g\tilde x)-\mathcal W_v(g\tilde y)\|^2)^{\frac12}\leq
\liminf_k(\sum_{g\in G}\|M^kW_{s_0}s(g\tilde x)-M^kW_{s_0}s(g\tilde y)\|^2)^{\frac12}$$
$$\leq\mbox{const } (\|s\|_\infty d(\tilde x,\tilde y)+\|s(p\tilde x)-s(p\tilde y)\|).$$
This shows that $p\tilde x\mapsto\ip{\mathcal W_vs}{\mathcal W_vs}'^{\frac12}$ is continuous, so 
$\ip{\mathcal W_vs}{\mathcal W_vs}'$ is, and this implies that $\mathcal W_vs\in\Xi$.
\par
The refinement equation is clearly satisfied because $\lim_kM^kW_{s_0}s(\tilde x)=(M\lim_k W_{s_0}s)(\tilde x)$.
\par
The last statement follows from the fact that $\mathcal P(x_0)$ is the orthogonal projection onto the eigenspace $E_1$.
\end{proof}

From the ``scaling functions'' $\mathcal W_v$ defined in Theorem \ref{prop5-1}, we can construct the harmonic functions $h_{v_1,v_2}$ as the correlations $\mathcal W_{v_2}^*\mathcal W_{v_1}$ between these scaling functions. 

\begin{theorem}\label{thh-v}
For $v_1,v_2\in\bc^n$ define $h_{v_1,v_2}:=\mathcal W_{v_2}^*\mathcal W_{v_1}$. Denote by $h_{v_1}:=h_{v_1,v_1}=\mathcal W_{v_1}^*\mathcal W_{v_1}$.
\begin{enumerate}
\item The function $h_{v_1,v_2}$ is continuous and harmonic (i.e., $Rh_{v_1,v_2}=h_{v_1,v_2}$). 
\item For all $v\in E_1$, with $\|v\|=1$, and all continuous harmonic functions $h\geq 0$, with $h(x_0)v=v$, one has $h_v\leq h$. 
\item For all continuous harmonic functions $h$, $h(x_0)$ commutes with the projection $\mathcal P(x_0)$ onto the eigenspace $E_1$ (see also Theorem \ref{prop5-1}). 
\item The map $\Psi_{x_0}: h\mapsto h(x_0)\mathcal P(x_0)$ is a surjective $*$-algebra morphism between the algebra of continuous harmonic functions and $M_l(\bc)$.
\item For all $v_1,v_2\in E_1$, one has $h_{v_1,v_2}(x_0)w=\ip{v_1}{w}v_2$ for all $w\in E_1$.
\item There exist $n_2,...,n_p\geq 1$ and a $C^*$-algebras isomorphism $\alpha:\mathfrak{H}\rightarrow M_l(\bc)\oplus M_{n_2}(\bc)\oplus...\oplus M_{n_p}(\bc)$, such that $\mbox{proj}_1\circ\alpha(h)=\Psi_{x_0}(h)=h(x_0)\mathcal P(x_0)$ for all $h\in\mathfrak{H}$.
\item If $v\in E_1$ with $\|v\|=1$ then $h_v$ is a minimal projection in $\mathfrak H$. If $v_1,v_2\in E_1$ with $\|v_1\|=\|v_2\|=1$ and $v_1\perp v_2$ then $h_{v_1}*h_{v_2}=0$ in $\mathfrak H$.
\item If $v_1,v_2\in E_1$ then $\tau_{v_1,v_2}(f)=\ip{v_2}{f(x_0)v_1}$ is a linear functional on $End(\xi)$ which is invariant for the transfer operator, i.e., $\tau_{v_1,v_2}(Rf)=\tau_{v_1,v_2}(f)$ for all $f\in End(\xi)$. If $\{v_1,...,v_l\}$ is an orthonormal basis for $E_1$ then $\tau_{v_i,v_j}(h_{v_{i'},v_{j'}})=\delta_{ii'}\delta_{jj'}$.\
Moreover $\tau_{ii}$ is a pure state on $\mathfrak H$ for all $i\in\{1,...,l\}$.
\end{enumerate}

\end{theorem}

\begin{proof}
(i) The continuity of $h_{v_1,v_2}$ follows from the fact that $\mathcal W_{v_1}$ and $\mathcal W_{v_2}$ are in $\Hom{S,\Xi}$. Using the intertwining relation in Theorem \ref{thinter}, we have 
$$R(h_{v_1,v_2})=R(\mathcal W_{v_2}^*\mathcal W_{v_1})=(M\mathcal W_{v_2})^*M\mathcal W_{v_1}=\mathcal W_{v_2}^*\mathcal W_{v_1}=h_{v_1,v_2}.$$
\par
(ii) Since $h\geq 0$, we have $h^{1/2}(x_0)v=v$ and $h^{1/2}\in\End\xi$. Then, as in the proof of Theorem \ref{prop5-1}, taking $s_0\in S_1$ with $s_0(x_0)=v$, and $\|s_0\|_\infty\leq 1$, one has that for all $s\in S$,
 $\{M^k(W_{s_0}h^{1/2})s\}_k$ converges uniformly on compact sets to $\mathcal W_v$. Then, using Fatou's lemma:
 $$\ip{s}{\mathcal{W}_v^*\mathcal W_vs}=\ip{\mathcal W_vs}{\mathcal W_vs}'\leq\liminf_k\ip{M^k(W_{s_0}h^{1/2})s}{M^k(W_{s_0}h^{1/2})s}'=\liminf_k\ip{s}{R^k(h^{1/2}W_{s_0}^*W_{s_0}h^{1/2})s}$$$$\leq \liminf_k\ip{s}{R^k(h^{1/2}\|W_{s_0}^*W_{s_0}\|h^{1/2})s}\leq\ip{s}{hs},$$
 because $$\ip{W_{s_0}s}{W_{s_0}s}'=\ip{f}{f}'|\ip{s_0}{s}|^2\leq\ip{f}{f}'\|s_0\|_\infty^2\|s\|_\infty^2\leq\|s\|_\infty^2$$
 so $\|W_{s_0}^*W_{s_0}\|\leq 1$.
 \par
 This implies that $h_v\leq h$. 
 \par
 (iii) Let $h$ be a continuous harmonic function. Using Lemma \ref{lemmgx-0}, we have for all $v\in E_1$, $k\geq1$:
 $$h(x_0)v=R^kh(x_0)v={m^{(k)}}^*(\psi_1^k\tilde x_0)h(p\psi_1^k\tilde x_0)m^{(k)}(\psi_1^k\tilde x_0)v=(m(x_0)^k)^*h(x_0)m(x_0)^kv=(m(x_0)^k)^*h(x_0)v$$
Letting $k\rightarrow\infty$, and using Theorem \ref{prop5-1}, we obtain
\\$h(x_0)v=\mathcal P(x_0)^*h(x_0)v$, and this implies $h(x_0)\mathcal P(x_0)=\mathcal P(x_0)h(x_0)\mathcal P(x_0)$. 
Apply this equality to the harmonic function $h^*$, and take the adjoint to obtain $\mathcal P(x_0)h(x_0)=\mathcal P(x_0)h(x_0)\mathcal P(x_0)$. Thus $h(x_0)$ commutes with $\mathcal P(x_0)$.
\par
(iv) The map is clearly linear and preserves the adjoint. We prove that it preserves multiplication: take two continuous harmonic functions $a$, $b$ and $v\in E_1$. Then we have $a(x_0)b(x_0)v\in E_1$. Using a similar argument as above, we obtain
$$a*b(x_0)v=\lim_{k\rightarrow\infty}\frac1k\sum_{j=1}^kR^j(ab)(x_0)v=\lim_k\frac1k\sum_{j=1}^k(m(x_0)^j)^*a(x_0)b(x_0)v=a(x_0)b(x_0)v,$$
and we used Lemma \ref{lemeigmx-0} for the last equality. This shows that $a*b(x_0)\mathcal P(x_0)=a(x_0)b(x_0)\mathcal P(x_0)=a(x_0)\mathcal P(x_0) b(x_0)\mathcal P(x_0)$. Therefore $\Psi_{x_0}$ is $*$-algebra morphism.
\par
To see that $\Psi_{x_0}$ is surjective, let $b\in M_l(\bc)$. Pick a continuous $a\in \End\xi$ such that $a(x_0)v=bv$ for all $v\in E_1$. Then for $v\in E_1$:
$$T_1(a)(x_0)v=\lim_k\frac1k\sum_{j=1}^k R^j(a)(x_0)v=\lim_k\frac1k\sum_{j=1}^k(m(x_0)^j)^*a(x_0)v=bv$$
Therefore $\Psi(x_0)(T_1(a))=b$, and $\Psi_{x_0}$ is surjective.
\par
(v) Let $s_1,s_2\in S$ with $s_1(x_0)=w_1\in E_1$, $s_2(x_0)=w_2\in E_1$. Then
$$\ip{w_1}{h_{v_1,v_2}(x_0)w_2}=\ip{\mathcal W_{v_2}s_1}{\mathcal W_{v_1}s_2}'(x_0)=\sum_{g\in G}\ip{\mathcal P(g\tilde x_0)w_1}{v_2}\ip{v_1}{\mathcal P(g\tilde x_0)w_2}\ip{f(\tilde x_0)}{f(\tilde x_0)}$$
and using (\ref{eqpgx-0}) and Theorem \ref{prop5-1},
$$=\ip{w_1}{v_2}\ip{v_1}{w_2}=\ip{w_1}{\ip{v_1}{w_2}v_2}.$$
This proves (v).
\par
(vi) Since $\mathfrak{H}$ is a finite dimensional $C^*$-algebra, there exist $n_1,...,n_p$ and an isomorphism $\beta:\mathfrak H\rightarrow\oplus_{i=1}^pM_{n_i}(\bc)$. Consider the kernel $K:=Ker(\Psi_{x_0}\circ\beta^{-1})$. This is an ideal in $\oplus_iM_{n_i}(\bc)$. It is easy to see that $K=\oplus_i K_i$, where $K_i$ is an ideal in $M_{n_i}(\bc)$. Then $K_i=0$ or $K_i=M_{n_i}(\bc)$. Denote the complements of these ideals $K_i^c:=M_{n_i}(\bc)\ominus K_i$. Then $\Psi_{x_0}\circ\beta^{-1}$ restricts to an isomorphism on $\oplus_iK_i^c$. Since $M_l(\bc)$ has trivial center, only one of the summands $K_i^c$ is non-trivial. After somer relabeling we may assume $K_1^c$ is the non-trivial one, so $n_1=l$. Therefore  $\gamma_1:x_1\mapsto \Psi_{x_0}\circ\beta^{-1}(x_1,0,\cdots,0)$ is an isomorphism of $M_l(\bc)$ onto itself. Then, we have that $$\Psi_{x_0}\circ\beta^{-1}(\gamma_1^{-1}x_1,x_2,\cdots,x_p)=\Psi_{x_0}\circ\beta^{-1}(\gamma_1^{-1}x_1,0,\cdots,0)=\gamma_1(\gamma_1^{-1}x_1)=x_1$$ for all $(x_1,...,x_p)\in\oplus_iM_{n_i}$. Then, define $\alpha(h)=(\gamma_1\beta_1(h),\beta_2(h),\cdots,\beta_p(h))$ for $h\in\mathfrak H$. Then $\alpha$ is an isomorphism, and if $\alpha(h)=:(x_1,...,x_p)$ then $\beta_1(h)=\gamma_1^{-1}x_1$, and $\beta_2(h)=x_2,\cdots,\beta_p(h)=x_p$, so $\Psi_{x_0}h=\Psi_{x_0}\beta^{-1}(\gamma_1^{-1}x_1,x_2,\cdots,x_p)=x_1=\alpha_1(h)$. This gives the desired isomorphism $\alpha$.
\par
(vii) Consider $\alpha(h_v)$. If $b=(b_1,...,b_p)\in M_l(\bc)\oplus M_{n_1}(\bc)...\oplus M_{n_p}(\bc)$ is positive, and $b_1v=v$, then $\alpha^{-1}(b)$ is positive, and $\alpha^{-1}(b)(x_0)v=b_1v=v$. Using (ii), we obtain that $\alpha^{-1}(b)\geq h_v$, so $b\geq \alpha(h_v)$. Since we also have $\mbox{proj}_1\alpha(h_v)v=v$, we obtain that $\alpha(h_v)$ is a minimal projection. Therefore $h_v$ is. 
\par
Since $\alpha(h_{v_1})$ and $\alpha(h_{v_2})$ are minimal projections onto $v_1$ and $v_2$ respectively, the affirmation follows. 
\par
(viii) We have, with Lemma \ref{lemmgx-0}, for all $g\in\End\xi$,
$$\tau_{v_1,v_2}(Rg)=\sum_{j=1}^k\ip{v_2}{m^*(\psi_jx_0)g(\psi_jx_0)m(\psi_jx_0)v_1}=$$$$\ip{ m(x_0)v_2}{g(x_0) m(x_0)v_1}=\ip{v_2}{g(x_0)v_1}.$$
The other statements follows from (v) and (vi).
\end{proof}

\subsection{Strong convergence of the cascade algorithm}
 We show here that, when the peripheral spectrum of $R$ is completely described by Theorem \ref{thh-v}, then the iterates of the refinement operator converge strongly to the fixed points $\mathcal{W}_v$. We begin with a result that applies to a more general situation, and gives a sufficient condition for strong convergence.

\begin{theorem}\label{thstrconv}
Assume that
 $1$ is the only eigenvalue of $R$ on the unit circle. Let $W \in \text{Hom}_{C(X)}(S,\Xi)$, $Ws=\sum_{j=1}^k\ip{s_j}{s}f_j$ with $s_j\in S_1$, $f_j\in\Xi_1$ for all $j\in\{1,...,k\}$, and let
  $\tau_1, \dots, \tau_l$ be a family of continuous linear functionals on $\End\xi$ invariant under $R$ and that separates points in $\mathfrak H$. 
  If 
  $\tau_i((M^kW-W)^*(M^kW -W)) = 0$ for every $k \geq 0$ and every $i$ then 
  $M^k W$ converges with respect to the operator norm on $\text{Hom}_{C(X)}(S, \Xi)$
  and 
  $$ (\lim_kM^k W)^*(\lim_kM^kW) = T_1(W^*W).$$
   \begin{proof}
    The intertwinning relation  gives us 
    $$
    \|\langle  M^{k+l} W s - M^l W s,  M^{k+l} W s - M^l W s \rangle'  \|_\infty\\
      =  \| \langle ( M^k -1)M^l Ws, ( M^k -1)M^l Ws \rangle' \|_\infty$$
    $$=  \| \langle s , R^l( (M^kW -W)^*(M^k W - W)) s \rangle \|_\infty 
    \leq  \| \langle s, s \rangle \|_\infty \| R^l(  (M^kW -W)^*(M^k W - W))\|_\infty.$$
    
    Moerover, since 
     $\tau_i((M^kW-W)^*(M^kW -W)) = 0$ for every $k \geq 0$ and $\tau_i$ is invariant for $R$, we get that $\tau_i(T_1( (M^kW-W)^*(M^kW -W)) = 0$, and since $\tau_i$ are separating, we obtain $T_1((M^kW-W)^*(M^kW-W))=0$. 
    Since $R$ is quasicompact with no other eigenvalues on the unit circle than $1$,
    we obtain a decomposition $R = T_1 + S$ as operators on $L$, so $ST_1=T_1S=0$ and $R^n=T_1+S^n$ for all $n$. Since $T_1$ is compact, the essential spectral radius of $S$ is the same as the one of $R$ (see \cite{Nussbaum}), but as $S$ does not have eigenvectors for the eigenvalue $1$, the spectrum of $S$ is contained 
    in a disc of radius $1-\epsilon$ with  $\epsilon > 0$. By the spectral radius formula, we see that 
    $\lim_n \| S^n \|_L = 0$.
    Since $ \sup_k \|(M^kW-W)^*(M^kW -W)\|_1 < \infty$ (see Lemma \ref{lemlipmk}), we see that 
    $$
      \| R^n(  (M^kW -W)^*(M^k W - W))\| 
      \leq     \| R^n(  (M^kW -W)^*(M^k W - W))\|_L $$$$
     =  \| S^n(  (M^kW -W)^*(M^k W - W))\|_L \leq C \| S^n \|_L, $$
   
     for every $k$,   i.e. $M^kW$ form a Cauchy sequence in $\text{Hom}_{C(X)}(S, \Xi)$.
     
     The last claim is true since
$$
       \langle \lim_k M^kW s_1, \lim_k M^k W s_2 \rangle'  = \lim_k \langle M^kW s_1, M^k W s_2 \rangle' 
        = \lim_k \langle s_1, R^k (W^*W ) s_2 \rangle $$$$
        = \lim_k \langle s_1, T_1(W^*W)s_2 + S^k (W^*W) s_2 \rangle 
        = \langle s_1, T_1(W^*W)s_2 \rangle.$$

  \end{proof}

\end{theorem}

\begin{theorem}\label{computation}
Assume $R1=1$, $m(x_0)$ satisfies the $E(l)$ condition, $\dim\mathfrak H=l^2$ and $1$ is the only eigenvalue for $R$ on the unit circle. Let $\{v_1,\cdots,v_l\}$ be an orthogonal basis for $E_1=\{v\in\bc^d\,|\, m(x_0)v=v\}$ and let $s_1,\cdots,s_l\in S_1$ with $s_j(x_0)=v_j$ for all $1\leq j\leq l$. Define $W_{s_1},\cdots,W_{s_l}$ as in Theorem \ref{prop5-1}

\begin{enumerate}
\item
For all $1\leq j\leq l$, $M^kW_{s_j}$ converges to $\mathcal W_{v_j}$ with respect to the norm in $\Hom{S,\Xi}$;
\item 
For all $1\leq j_1,j_2\leq l$, 
$$h_{v_{j_1},v_{j_2}}:=\mathcal W_{v_{j_2}}^*\mathcal W_{v_{j_1}}=T_1(W_{s_{j_2}}^*W_{s_{j_1}}).$$
\item
The algebra $\mathfrak H$ is isomorphic to $M_l(\bc)$, and 
$$h_{v_1}+\cdots...+h_{v_l}=1.$$
\end{enumerate}
\end{theorem}
\begin{proof}
  We can use Theorem \ref{prop5-1} and we have that $M^kW_{s_j}$ converges pointwise to $\mathcal W_{v_j}$. Also, $h_{v_i,v_j}$ form a basis for $\mathfrak H$ and $\tau_{i,j}$ form a family of functionals invariant under $R$ that separate the points in $\mathfrak H$. \par We use Theorem \ref{thstrconv}, and we perform the following computation:
  since $s_i(x_0)=v_i$, we have:
\begin{align*}
& \tau_{i,i} (  (M^{k} W_{s_j}  - W_{s_j})^* (M^k W_{s_j} - W_{s_j}) ) 
  =  \langle (M^k W_{s_j} - W_{s_j}) s_i , (M^k W_{s_j} - W_{s_j}) s_i \rangle' (x_0) \\
  = &  \sum_{g \in G} | \langle s_j \circ p ( \tilde r ^{-k} g \tilde x_0) , 
  m^{(k)}(\tilde r ^{-k}  g x_0) s_i (x_0) \rangle f( \tilde r ^{-k} g \tilde x_0) 
   - \langle s_j(x_0), s_i(x_0) \rangle f(g\tilde x_0) |^2
\end{align*}
(see also \eqref{eqmkw}).\par
For $g\in G$, we have a unique factorization 
 $$
 \tilde r ^{-k} g = g' \psi_{\omega_k} \dots \psi_{\omega_1}, 
$$
 where $g' \in G$. The elements $g \in A^k G$ are exactly those for which the above factorization is 
 $g'    \psi_1 \dots \psi_1  $.
 Moreover, with Lemma \ref{lemmgx-0}, we have that 
 $m^{(k)}(  \psi_{\omega_k} \dots \psi_{\omega_1} \tilde x_0) v_i = \delta_{ \omega ,1\cdots 1} v_i $. Also $\|f(\tilde x_0)\|=1$ and $f(g\tilde x_0)=0$ for $g\neq 0$.
This implies that the above expression equals
 \begin{align*}
   &  \sum_{g' \in G} \sum_{\omega \in \Omega_k} 
   | \langle s_j \circ p ( \psi_{\omega_k} \dots \psi_{\omega_1} \tilde x_0) , 
  m^{(k)}(   \psi_{\omega_k} \dots \psi_{\omega_1} \tilde x_0) s_i (x_0) \rangle 
f( g'  \psi_{\omega_k} \dots \psi_{\omega_1} \tilde x_0) 
  \\ & - \langle s_j(x_0), s_i(x_0) \rangle f(\tilde r^{k}g'  \psi_{\omega_k} \dots \psi_{\omega_1}   \tilde x_0) |^2 \\
    = &    \sum_{g \in G} \sum_{\omega \in \Omega_k} |\langle s_j(x_0) ,  s_i (x_0) \rangle |^2
    |  
  \delta_{1\cdots 1, \omega} f( g  \psi_{\omega_k} \dots \psi_{\omega_1} \tilde x_0) 
   -  f(\tilde r^kg' \psi_{\omega_k} \dots \psi_{\omega_1}   \tilde x_0) |^2 \\
    = &    \delta_{i,j}|f(x_0)-f(x_0)|=0,\quad (\mbox{since }f(gx_0)=0\mbox{ for }g\neq1). 
 \end{align*}
 Then for $i\neq i'$, by an application of Cauchy-Schwarz, we have:
 \begin{align*}
   & |\tau_{i,i'}( ( M^k W_{s_j} - W_{s_j})^*(M^k W_{s_j}- W_{s_j}) ) |^2 \\ & \leq   \tau_{i,i}( ( M^k W_{s_j} - W_{s_j})^*(M^k W_{s_j}- W_{s_j}) )  \tau_{j,j}( ( M^k W_{s_j} - W_{s_j})^*(M^k W_{s_j}- W_{s_j}) )=0.
 \end{align*}
 Thus $M^kW_{s_j}$ converges strongly, by Theorem \ref{thstrconv} and (i) follows.
\par
(ii) follows directly from (i) as in Theorem \ref{thstrconv}.
\par
(iii) follows directly from Theorem \ref{thh-v}.
\end{proof}
\begin{acknowledgements}The authors are pleased to acknowledge helpful
discussions with Professors Ola Bratteli, Deguang Han, Palle Jorgensen, Sergey Neshveyev, Qiyu Sun.
\end{acknowledgements}

\bibliography{dreftr}

\end{document}